\documentclass[11pt]{amsart}
\usepackage{amssymb}
\usepackage{amsmath, amsthm, amsopn, amsfonts}
\usepackage[dvips]{graphicx}
\usepackage{graphpap}
\usepackage[francais, english]{babel}
\def\R{\mathbb R}

\def\T{\mathbb T}
\def\N{\mathbb N}

\newtheorem{Theorem}{Theorem} 

\newtheorem{lemma}{Lemma} 
\newtheorem{prop}{Proposition}

\newcommand{\rien}[1]{\relax}

\newcommand{\ep}{\varepsilon} 
\newcommand{\tu}{u^\varepsilon} 
\renewcommand{\u}{u^\varepsilon}
\newcommand{\vp}{\varphi^\varepsilon}
\renewcommand{\ep}{\varepsilon} 
\newcommand{\intR}{\int_{\R^2}}
\newcommand{\p}{\partial} 
\newcommand{\pa}{\partial^{-1}_x}

\newcommand{\re}[1]{(\ref{#1})} 
 \newcommand{\lab}[1]{\label{#1}}

\newcommand{\intr}{\int}
 
\newcommand{\beq}{\begin{equation}}
\newcommand{\eneq}{\end{equation}} 
\newcommand{\beqn}{\begin{eqnarray}}
\newcommand{\eneqn}{\end{eqnarray}} 
\newcommand{\beqnon}{\begin{eqnarray*}}
\newcommand{\eneqnon}{\end{eqnarray*}} 
\newcommand{\tendsto}[1]{ \renewcommand{\arraystretch}{0.5}
\begin{array}[t]{c}
 \longrightarrow \\
 { \scriptstyle #1 }
\end{array}
\renewcommand{\arraystretch}{1}}

\begin{document}

\selectlanguage{english}

\title
[ KP-I on the background of a non localized solution ]
{\bf Global well-posedness for the KP-I equation on the background of a non localized solution }
\author{L.~Molinet} 
\address{L.A.G.A., Institut Galil\'ee, Universit\'e Paris 13, 93430 Villetaneuse }
\author{J. C.~Saut}
\address{ Universit\'e de Paris-Sud, UMR de Math\'ematiques, B\^at. 425, 91405 Orsay Cedex }
\author {N.~Tzvetkov}
\address
{D\'epartement de Math\'ematiques, Universit\'e Lille I, 59 655 Villeneuve d'Ascq Cedex}
\begin{abstract}
We prove that the Cauchy problem for the KP-I equation is globally well-posed
for initial data which are localized perturbations (of arbitrary size) of a
{\bf non-localized} (i.e. not decaying in all directions) traveling wave
solution (e.g. the KdV line solitary wave or the Zaitsev solitary waves which
 are localized in $x$ and $y$ periodic or conversely).
\end{abstract}  
\maketitle
\section{Introduction}
We study here the initial value problem for the Kadomtsev-Petviashvili (KP-I) equation
\begin{equation}
(u_t+u_{xxx} +u u_x )_x - u_{yy} =0, 
\label{1.1}
\end{equation}
where $u=u(t,x,y)$, $(x,y)\in\R^2$, $t\in\R$, with initial data
\begin{equation}\label{1.2}
u(0,x,y)=\phi(x,y)+\psi_c(x,y),
\end{equation}
where $\psi_c$ is the profile\footnote{This means that $\psi(x-ct,y)$ solves
(\ref{1.1}).} of a non-localized (i.e. not decaying in all spatial directions)
traveling wave of the KP-I equation moving with speed $c\neq 0$.
This $\psi_c$ could be for instance the line soliton of
the Korteweg- de Vries (KdV) equation
\begin{equation}\label{KdV}
\psi_c(x,y)=3c\,{\rm cosh}^{-2}\Big(\frac{\sqrt{c}\, x}{2}\Big)
\end{equation}
In (\ref{KdV}) the KdV soliton is of course considered as a two dimensional (constant
in $y$) object. 
Another possibility to see (\ref{KdV}) as a solution of KP-I is to consider
(\ref{1.1}) posed on $\R\times\T$. Global solutions of (\ref{1.1}) for data on
$\R\times\T$, including data close to (\ref{KdV}) were recently constructed in a work by 
Ionescu-Kenig~\cite{IK}.
In (\ref{1.2}), the function $\psi_c$ may also be the profile of the Zaitsev \cite{Z} traveling waves (see also \cite{TM})
which is localized in $x$ and periodic in $y$ :
\begin{equation}\label{1.3}
\psi_c(x,y)=12\alpha^{2} 
\frac{1-\beta\cosh(\alpha x)\cos(\delta y)}
{(\cosh(\alpha x)-\beta\cos(\delta y))^2},
\end{equation}
where
$$
(\alpha,\beta)\in
]0,\infty[\times ]-1,1[,
$$
and the propagation speed is given by
$$
c=\alpha^{2}\frac{4-\beta^2}{1-\beta^2}\,.
$$
Let us observe that the transform $\alpha\rightarrow i\alpha$,
$\delta\rightarrow i\delta$, $c\rightarrow ic$ produces solutions of
(\ref{1.1}) which are periodic in $x$ and localized in $y$. The profiles of
these solutions are also admissible in (\ref{1.2}), under the assumption
$|\beta|>1$. Notice that for $\beta=0$, (\ref{1.3}) coincides with (\ref{KdV}).
\\

The global well-posedness of (\ref{1.1})-(\ref{1.2}) with data given by
(\ref{KdV}) which will be proved in this paper can be viewed as a preliminary step towards the rigorous
mathematical justification of the (conjectured) nonlinear instability of the KdV soliton
with respect to transversal perturbations governed by the KP-I flow. This
question is, as far as we know, still an open problem (see however \cite{APS}
for a linear analysis of the instability and \cite{GHS} for a linear
instability analysis in the framework of the full Euler system).
The instability scenario of the line soliton seems to be a symmetry breaking
phenomenon :
the line soliton should evolve towards the Zaitsev solitary wave
(\ref{1.3}). Note that Haragus and Pego \cite{HP} have shown that this
solution is the only one close to the
line soliton which is periodic in $y$ and decays to zero as $x\rightarrow
\infty$.
\\

The question of solving (\ref{1.1}) together with the initial data (\ref{1.2})
when $\psi$
is the profile of the KdV line soliton, has been recently addressed 
by Fokas and Pogrobkov \cite{FP}, by the inverse scattering transform (IST)
techniques.
However, the Cauchy problem is not rigorously solved in \cite{FP} and it is
unlikely that it could be solved for an arbitrary large data $\phi$ using
IST since the Cauchy problem with purely localized data has been solved by IST
techniques only for {\bf small} initial data (see \cite{W,Zh}).
\\

On the other hand, PDE techniques have been recently fruitfully used to obtain
the global solvability of the KP equation with arbitrary large initial data,
starting with the pioneering paper of Bourgain \cite{Bo} on the KP-II equation
(that is (\ref{1.1}) with $+u_{yy}$ instead of $-u_{yy}$). 
In \cite{Bo}, the global well-posedness of the KP-II equation for data in
$H^s(\R^2)$, $s\geq 0$ is established.
The result in \cite{Bo} is obtained by performing the Picard iteration scheme
to an equivalent integral equation in the Fourier transform restriction spaces
of Bourgain.
The situation for the KP-I equation turned out to be more delicate. 
We showed in \cite{MST-bis} that the Picard iteration scheme can not be
applied in the context of the KP-I equation as far as one considers initial
data in Sobolev spaces. Sobolev spaces are natural, since the
conservation laws for the KP-I equation control Sobolev type norms.
In \cite{KT1}, a quite flexible method is introduced that allows to
incorporate the dispersive effects in a context of a compactness method for proving
the well-posedness. The work \cite{KT1} was in turn inspired by the
considerations in \cite{BGT} in the context of the Nonlinear Schr\"odinger
equation on a compact manifold which is solved in \cite{BGT} as a semi-linear
problem (i.e. by the Picard iteration scheme). The main point in \cite{KT1} is to realize that the idea of
\cite{BGT} can also be used in the context of a quasi-linear problem.
The method of \cite{KT1} turned out to be useful in the context of the KP-I
equation (and some other models such as the Schr\"odinger maps \cite{JKato,KN}) and the first global well-posedness result
for the KP-I equation has been obtained by the authors of the present paper in \cite{MST}. 
This result has been improved (i.e. the space of the allowed initial data is larger) by Kenig~\cite{K}. 
Together with the idea of \cite{KT1}, a new commutator estimate for the KP-I
equation is used in \cite{K}. The main point in
Kenig's result is the proof that KP-I is locally well-posed for data in the 
space
$$
\{
u\in L^2(\R^2)\, :\, \partial_x^{-1}u_y\in L^2(\R^2),\quad |D_x|^s u\in
L^2(\R^2), \,\, s>3/2 
\}
$$
All papers \cite{Bo}, \cite{MST} and \cite{K} consider the KP equations in spaces of ``localized'' (zero
at infinity) functions.
\\

The main goal of this paper is  
to prove that, for a large class of $\psi_c$,  the Cauchy problem (\ref{1.1}), (\ref{1.2}) is globally
well-posed for data $\phi\in Z$, where
$$
Z:= \{ u\in L^2(\R^2)\, :\, \partial_x^{-2}u_{yy}\in L^2(\R^2),\quad u_{xx}\in
L^2(\R^2)
\}.
$$
(notice that $u\in Z$ implies $u_y\in L^2(\R^2)$ and $\partial_{x}^{-1}u_{y}\in
L^2(\R^2)$).
More precisely, we have the following result.
\begin{Theorem}\label{thm1}
Let $\psi_c(x-ct,y)$ be a solution of the KP-I equation such that 
$$
\psi_c\,:\, \R^2\longrightarrow \R
$$ 
is bounded with all its derivatives\footnote{The bounds can of course depend on the propagation speed $c$.}.
Then for every $\phi\in Z$ there exists a unique  global solution
$u$ of \re{1.1} with initial data \re{1.2} satisfying for all $T>0$,
$$
[u(t,x,y)-\psi_c(x-ct,y)]\in C([0,T];Z),\quad  \partial_x
[u(t,x,y)-\psi_c(x-ct,y)]\in L^1_T L^\infty_{xy}\, . 
$$
Furthermore, for all $ T> 0$,  the map $ \phi \mapsto u $ is continuous
from $Z$ to $ C([0,T];Z))$.
\end{Theorem}
Since the KP-I equation is time reversible, a similar statement to
Theorem~\ref{thm1} holds for negative times as well.
Let us notice that the assumptions on $\psi_c$ in Theorem~\ref{thm1} are clearly satisfied by the
line or the Zaitsev solitary wave.
In the proof of Theorem~\ref{thm1}, we write the solution $u$ of (\ref{1.1}), (\ref{1.2}) as
$$
u(t,x,y)=\psi_c(x-ct,y)+v(t,x,y)
$$
where $v$ is localized. This $v$ satisfies the equation
 \begin{equation}
(v_t+v_{xxx} +v v_x +\partial_x (\psi_c v) )_x - v_{yy} =0,\quad v(0,x,y)=\phi(x,y). \label{LineKP}
\end{equation}
Our strategy is then to adapt the proof of \cite{MST,K}. Starting from the
local well-posedness result, we implement a compactness method based on
``almost conservation laws''. 
New terms occur with respect to \cite{MST} but they are controlled since $\psi$ and its derivatives are bounded.
It is of importance for our analysis that equation (\ref{LineKP}) does not
contain a source term.
\\

We refer to the work by Gallo~\cite{Gallo} and the references therein, where non vanishing
at infinity solutions to one dimensional dispersive models are constructed.
\\

Let us notice that the framework considered in Theorem~\ref{thm1} is also a
convenient one for a rigorous study of the interaction of a line and lump
solitary waves (see e.g. \cite{Fre}).
\\

The rest of this paper is organized as follows. In the next section, using a
compactness method, we prove a basic well-posedness result for (\ref{LineKP}).
In Section~3, inspired by the formal KP-I conservation laws, we provide bounds 
for some Sobolev type norms of the local solutions. These bounds are however not
sufficient to get global solutions. For that reason, in Section~4 we prove a
Strichartz type bound. This bound is then used in Section~5 to get a first
global well-posedness result. In Sections 6 and 7 we extend the well-posedness
to the class $Z$ introduced above. The last section is
devoted to the ``usual'' KP-I equation. We show how the estimates of Section~4
can be used to give a slight improvement of the Kenig well-posedness result \cite{K}.  
\section{Local well-posedness for data in $H^{s}_{-1}(\R^2)$, $s>2$} 
In this section, we prove a basic local well-posedness result for
\re{LineKP}. The proof follows a standard compactness method.
We however need to work in Sobolev spaces of integer indexes, in order
to make work the commutator estimates related to the term $\partial_x (\psi_c v)$.
\\

By $H^s(\R^2)$, $s\in\R$, we denote the classical Sobolev spaces.
The local existence result for \re{integrateLineKP} will be obtained in the spaces
$H^{s}_{-1}(\R^2)$ equipped with the norm 
$$
\|u\|_{H^{s}_{-1}(\R^2)}=
\|(1+|\xi|^{-1})\langle|\xi|+ |\eta|\rangle^{s}
\widehat{u}(\xi,\eta)\|_{L^{2}_{\xi\eta}},
$$
where $\langle\cdot\rangle=(1+|\cdot|^2)^{\frac{1}{2}} $ and $\widehat{u}$
denotes the Fourier transform of $u$. The spaces $H^{s}_{-1}(\R^2)$ are
adapted to the specific structure of the KP type equations. In the estimates it will  appear
 $ W^{r,\infty} $-norms of
 $ \psi_c $ defined for any integer $ r\ge 0 $ by
$$
\|\psi_c\|_{W^{r,\infty}} =\sum_{0\le |\alpha|=|(\alpha_1,\alpha_2)|\le r} \|\partial_x^{\alpha_1} \partial_y^{\alpha_2} \psi_c\|_{L^\infty_{xy}}\quad .
$$

Consider the ``integrated'' equation  (\ref{LineKP})
\begin{equation}
u_t+u_{xxx} -\partial_x^{-1} u_{yy} + u u_x +\partial_x (\psi u)=0,  
\label{integrateLineKP}
\end{equation}
with initial data
\begin{equation}
u(0,x,y)=\phi(x,y). 
\label{integrateLineKPbis}
\end{equation}
For the solutions we study in this section, (\ref{LineKP}) may be substituted
by (\ref{integrateLineKP}).
For conciseness we skip the $c$ of $\psi_c$ in (\ref{integrateLineKP}) and we
suppose that $c=1$ in the sequel. Of course, the case $c\neq 1$ can be treated
in exactly the same manner.
We have the following local well-posedness result for \re{integrateLineKP}.
\begin{prop}\label{compact}
Let $ s>2 $ be an integer.  
Then for every $\phi \in H^s_{-1}(\R^2) $ there exists $T\gtrsim (1+\|\phi\|_{H^s})^{-1}$
and a unique solution $u$ to \re{integrateLineKP} on the time interval $[0,T]$ satisfying
\begin{equation*}
u\in C([0,T];H^{s}(\R^2))\,, \quad u\in L^{\infty}([0,T];H^{s}_{-1}(\R^2))\, .
\end{equation*}
In addition, for $t\in [0,T]$,
\begin{equation}\label{key}
\|u(t,\cdot)\|_{H^s(\R^2)}\leq C\|\phi\|_{H^s(\R^2)}\,\,
\exp\Big(c\|\nabla_{x,y}u\|_{L^{1}_{T}L^{\infty}_{xy}}+
cT\|\psi\|_{W^{s,\infty}}\Big).
\end{equation}
Moreover if $\phi \in H^{\sigma}_{-1}(\R^2) $ where $\sigma>s$ is an integer then 
$$
u\in C([0,T];H_{-1}^{\sigma}(\R^2)).
$$
Finally, the map $ \phi \mapsto u_\phi $ is continuous from $H^s(\R^2) $ to $ C([0,T];H^{s}(\R^2))$.
\end{prop}
\begin{proof}[Proof of Proposition~\ref{compact}.]
The process is very classical (see \cite{IN} for a closely related result). For $ \varepsilon >0 $ we look at
the regularized equation 
\begin{equation}\label{reg}
u^\ep_t +\ep \Delta^2 u^\ep_t=-u_{xxx}^\ep +\partial_x^{-1} u_{yy}^\ep - u^\ep
u^\ep_x -\partial_x (\psi u^\ep)\,
\end{equation}
where $\Delta=\partial_x^2+\partial_y^2$ is the Laplace operator. 
Equation (\ref{reg}) with initial condition 
$$
\phi_{\ep}=(1-\sqrt{\ep}\Delta)^{-1}\phi
$$
can be rewritten under the form 
\begin{equation}\label{reg-bis}
u^{\ep}(t)=L^{\ep}(t)\phi_{\ep}-\int_{0}^{t}\, L^{\ep}(t-\tau)
(1+\ep\Delta^2)^{-1}
\big( u^\ep(\tau)
u^\ep_x(\tau) +\partial_x (\psi(\tau) u^\ep(\tau) \big)d\tau\,,
\end{equation}
where 
$$
L^{\ep}(t):=\exp\Big(-t(1+\ep\Delta^2)^{-1}(\partial_{x}^{3}-\partial_x^{-1}\partial_y^2)\Big)\, .
$$
Notice that, thanks to the regularization effect of $(1+\ep\Delta^2)^{-1}$, for $\ep\neq 0$, the map
$$
u\longrightarrow (1+\ep\Delta^2)^{-1}
\big( u\, u_x +\partial_x (\psi u \big)\big)
$$
is locally Lipschitz from $H^s(\R^2) $ to $ H^s(\R^2)$, provided
$s>2$. The operator $L^{\ep}(t)$ is clearly bounded on $H^s(\R^2) $ and therefore
by the Cauchy-Lipschitz-Picard theorem there is a unique local solution 
\begin{equation}\label{picard}
u^{\ep}\in C([0,T];H^s(\R^2))
\end{equation}
of (\ref{reg-bis}) with 
$$
T\gtrsim 
(1+\|\phi_{\ep}\|_{H^s})^{-1}\geq
(1+\|\phi\|_{H^s})^{-1}\,.
$$ 
Thanks to the perfect derivative structure of the integral term in
(\ref{reg-bis}), we also obtain that
\begin{equation*}
u^{\ep}\in C([0,T];H^s_{-1}(\R^2)).
\end{equation*}
Notice also that, thanks to the assumption $s>2$, $(1-\sqrt{\ep}\Delta)u^{\ep}$ belongs to $H^s(\R^2)$.
We next study the convergence of $u^{\ep}$ as $\ep\rightarrow 0$. For that
purpose, we establish a priori bounds, independent\footnote{ Notice that the bounds on
$\|u^{\ep}(t,\cdot)\|_{H^s}$ resulting from the Cauchy-Lipschitz theorem
applied to (\ref{reg}) are unfortunately very poor (depending on $\ep$).}
of $\ep$ on
$\|u^{\ep}(t,\cdot)\|_{H^s}$ on time intervals of size of order
$(1+\|\phi\|_{H^s})^{-1}$. 
\\

Multiplying (\ref{reg}) with $u^{\ep}$, after an
integration by parts, we obtain that
\begin{equation}\label{101}
\frac{d}{dt}\Big[\|u^{\ep}(t,\cdot)\|_{L^2}^{2}+
\ep\|\Delta u^{\ep}(t,\cdot)\|_{L^2}^{2}\Big]\lesssim\|\psi_x\|_{L^{\infty}}\|u^{\ep}(t,\cdot)\|_{L^2}^{2}\, .
\end{equation}
Let us recall a classical commutator estimate (see e.g. \cite{KP}).
\begin{lemma} \label{comm}
Let $\Delta$ be the Laplace operator on $\R^n$, $n\geq 1$.
Denote by $J^{s}$ the operator $(1-\Delta)^{s/2}$.
Then for every  $s> 0$,
\begin{equation*}
\|[J^s,f] g\|_{L^2(\R^n)}
\lesssim  \|\nabla f\|_{L^\infty(\R^n)} \|J^{s-1} g\|_{L^2(\R^n)}+
\|J^s f\|_{L^2(\R^n)} \|g\|_{L^\infty(\R^n)} \, .
\end{equation*}
\end{lemma}
Using Lemma \ref{comm}, we obtain  for fixed $ y $,
$$
\|[\partial_x^s,u(x,\cdot)]u_x(x,\cdot)\|_{L^2_x}
\lesssim
\|u_x(x,\cdot)\|_{L^{\infty}_x}
\|J^{s}_x u(x,\cdot)\|_{L^{2}_x} \;
$$
 where $J_x^s=(1-\partial_x^2)^{s/2}$. Squaring, integrating over $ y $ and integrating by parts, it yields
\begin{equation}\label{kando1}
\Big|\intR \partial^s_x ( u u_x  ) \partial^s_x u \Big| \lesssim
\|u_x\|_{L^\infty_{xy}} \|J^s_x u\|_{L^2_{xy}}^2 \quad .
\end{equation}
 On the other hand, by Leibniz rule and integration by parts (recall that
 $ s$ is an integer), we get
\begin{equation}\label{kando2}
\Big|\intR \partial_x^s (\psi u )_x \partial^s_x u \Big|
\lesssim
\|\psi_x\|_{W^{s-1,\infty}}\|J^s_x u\|_{L^2_{xy}}^2 \quad .
\end{equation}
Therefore applying $\partial^s_x$ to 
(\ref{reg}), multiplying it with $\partial ^s_x u^{\ep}$  gives
\begin{multline}\label{103}
\frac{d}{dt}
\Big[
\|\partial_x^s u^{\ep}(t,\cdot)\|_{L^2}^{2}
+\ep
\|\partial_x^s \Delta u^{\ep}(t,\cdot)\|_{L^2}^{2}
\Big]
\lesssim
\\
\lesssim
\Big(\|u^{\ep}_x(t,\cdot)\|_{L^{\infty}}+\|\psi_x\|_{W^{s-1,\infty}}\Big)
\|J^s_x u^{\ep}(t,\cdot)\|_{L^2(\R^2)}^{2}
\end{multline}
Next we estimate the $y$ derivatives.
In the same way, using Lemma~\ref{comm}, we obtain that for a fixed $x$
\begin{multline}\label{y}
\|[\partial_y^s,u(x,\cdot)]u_x(x,\cdot)\|_{L^2_y}\lesssim
\\
\lesssim
\Big(\|u_x(x,\cdot)\|_{L^{\infty}_y}+\|u_y(x,\cdot)\|_{L^{\infty}_y}\Big)
\Big(\|J^{s-1}_y u_x(x,\cdot)\|_{L^{2}_y}+\|J^{s}_y u(x,\cdot)\|_{L^{2}_y}\Big),
\end{multline}
where $J^s_y=(1-\partial_y^2)^{s/2}$. 
Squaring (\ref{y}), integration over $x$ and an integration by parts yield
\begin{equation*}
\Big|\intR \partial^s_y ( u u_x  ) \partial^s_y u \Big| \lesssim
\|\nabla_{x,y}u\|_{L^\infty_{xy}} \|u\|_{H^s(\R^2)}^2 \quad .
\end{equation*}
 On the other hand, by Leibniz rule and integration by parts,
\begin{equation*}
\Big|\intR \partial_y^s (\psi u )_x \partial^s_y u \Big| 
\lesssim
\|\psi\|_{W^{s,\infty}}\|u\|_{H^s(\R^2)}^2 \quad .
\end{equation*}
Therefore applying $\partial_y^s$ to (\ref{reg}) and multiplying it by
$\partial_{y}^{s}u^{\ep}$ gives
\begin{multline}\label{104}
\frac{d}{dt}
\Big[
\|\partial_y^s u^{\ep}(t,\cdot)\|_{L^2}^{2}
+\ep
\|\partial_y^s\Delta u^{\ep}(t,\cdot)\|_{L^2}^{2}
\Big]
\lesssim
\\
\lesssim
\Big(\|\nabla_{x,y} u^{\ep}(t,\cdot)\|_{L^{\infty}}+\|\psi\|_{W^{s,\infty}}\Big)
\|u^{\ep}(t,\cdot)\|_{H^s(\R^2)}^{2}\,.
\end{multline}
Since for $s$ integer
$$
\|u\|_{H^s}\approx \|u\|_{L^2}+\|\partial^{s}_x u\|_{L^2}+\|\partial_y^s u\|_{L^2},
$$
combining (\ref{101}), (\ref{103}) and (\ref{104}) gives
\begin{multline*}
\frac{d}{dt}
\Big[
\|u^{\ep}(t,\cdot)\|_{H^s(\R^2)}^{2}
+\ep
\|\Delta u^{\ep}(t,\cdot)\|_{H^s(\R^2)}^{2}
\Big]
\lesssim
\\
\lesssim
\Big(\|\nabla_{x,y}u^{\ep}(t,\cdot)\|_{L^{\infty}}+ \|\psi\|_{W^{s,\infty}}\Big)
\|u^{\ep}(t,\cdot)\|_{H^s(\R^2)}^{2}
\end{multline*}
and therefore using that
$$
\|\phi_{\ep}\|_{H^s(\R^2)}^{2}+\ep\|\Delta\phi_{\ep}\|_{H^s(\R^2)}^{2}\leq
\|\phi\|_{H^s(\R^2)}^{2}
$$
by the Gronwall lemma for every $T>0$ on the time of existence
of $u^{\ep}$,
\begin{equation}\label{bez}
\|u^{\ep}\|_{L^{\infty}_{T}\,H^s(\R^2)}
\leq
\|\phi\|_{H^s(\R^2)}\exp\Big(c\|\nabla_{x,y} u^{\ep}\|_{L^{1}_{T}L^{\infty}_{xy}}+
cT\|\psi\|_{W^{s,\infty}}\Big)\, .
\end{equation}
which is the key inequality. Since $s>2$, using the Sobolev embedding, we get
\begin{equation}\label{sob}
\int_{0}^{T}\|\nabla_{x,y}u^{\ep}(\tau,\cdot)\|_{L^{\infty}}d\tau
\leq CT\|u^{\ep}(t,\cdot)\|_{L^{\infty}_{T}\, H^s(\R^2)}\, .
\end{equation}
Using (\ref{bez}), (\ref{sob}) and the continuity of $u^{\ep}(t)$ with
respect to time (see (\ref{picard})), we obtain that there exists
$C>0$ such that if
$$
T\lesssim (1+\|\phi\|_{H^s(\R^2)})^{-1}
$$
then
\begin{equation}\label{comp1}
\int_{0}^{T}\|\nabla_{x,y} u^{\ep}(\tau,\cdot)\|_{L^{\infty}(\R^2)}d\tau \leq C
\end{equation}
and
\begin{equation}\label{comp2}
\|u^{\ep}\|_{L^{\infty}_{T}H^{s}(\R^2)}\leq C\|\phi\|_{H^s(\R^2)}\, .
\end{equation}
We next estimate the anti-derivatives of $u^{\ep}$.
Let $v^{\ep}:=\partial_x^{-1}u^{\ep}$. Then, using (\ref{reg-bis}), we obtain
that $v^{\ep}$ solves the equation
\begin{equation*}
v^{\ep}(t)=L^{\ep}(t)\partial_{x}^{-1}\phi_{\ep}-
\int_{0}^{t}\, L^{\ep}(t-\tau)
(1+\ep\Delta^2)^{-1}
\big(\frac{1}{2}(u^\ep(\tau))^{2}+\psi(\tau) u^\ep(\tau) \big)d\tau\, .
\end{equation*}
Therefore, since $s>2$, using the Leibniz rule and the Sobolev inequality, we get the bound
\begin{equation}\label{pak}
\|v^{\ep}\|_{L^{\infty}_{T}H^s(\R^2)}\leq
\|\partial_{x}^{-1}\phi\|_{H^s(\R^2)}+
CT\|u^{\ep}\|_{L^{\infty}_{T}H^{s}(\R^2)}
\Big(\|u^{\ep}\|_{L^{\infty}_{T}H^{s}(\R^2)}+\|\psi\|_{W^{s,\infty}}\Big)\, .
\end{equation}
Coming back to the equation (\ref{reg}), we infer from (\ref{pak}) that the
sequence
$(\partial_{t}(u^{\ep}))$ is bounded in a weaker norm, say in $L^{\infty}_{T}
H^{s-2005}(\R^2)$. Therefore from the Aubin-Lions compactness theorem (see e.g. \cite{Lions}), we obtain that
$u^{\ep}$ converges, up to a subsequence, to some limit $u$ in the space $L^{2}_{loc}((0,T)\times \R^2)$
which satisfy (\ref{key}) and 
\begin{equation}\label{regu}
\int_{0}^{T}\|\nabla_{x,y} u(\tau,\cdot)\|_{L^{\infty}(\R^2)}d\tau \leq C\, .
\end{equation}
Thanks to (\ref{pak}), we obtain that (up to a subsequence)
$\partial_x^{-1}u^{\ep}$ converges in ${\mathcal D}'((0,T)\times \R^2)$
to a limit which can be identified as $\partial_{x}^{-1}u$.
By writing the nonlinearity $uu_x$ as $\frac{1}{2}\partial_x(u^2)$, 
passing into a limit in the equation (\ref{reg}) as $\ep\rightarrow 0$,
we obtain that the function $u$ satisfy the equation (\ref{integrateLineKP}) in
the distributional sense. Moreover thanks to (\ref{comp2}) and (\ref{pak}), we
obtain that
$$
u\in L^{\infty}([0,T];H^{s}_{-1}(\R^2))\, .
$$
In addition, thanks to the Lebesgue theorem, $\phi_{\ep}$ converges in $H^{s}_{-1}(\R^2)$ to $\phi$
as $\ep\rightarrow 0$ and thus $u$ satisfies the initial condition (\ref{integrateLineKPbis}).
Next, if $\phi\in H^{\sigma}_{-1}(\R^2)$ with $\sigma\geq s$, $\sigma\in \N$, then
as above, we get the estimate (\ref{key}) with $\sigma$ instead of $s$ which,
in view of (\ref{regu}), yields the propagation of the $H^{\sigma}(\R^2)$
regularity. We next estimate the anti-derivatives of $u$ in $H^{\sigma}$ by
invoking (\ref{pak}) (with $\sigma$ instead of $s$) in the limit $\ep\rightarrow 0$.
The uniqueness is straightforward from the
Gronwall lemma and (\ref{regu}).
The continuity of the flow map and the fact that the solution
is a continuous curve in $H^s(\R^2)$ can be obtained by the Bona-Smith
approximation argument \cite{BS}. We do not give the details of this
construction in this section since a completely analogous discussion will
be performed later in this paper. 
\end{proof}
Let us next state a corollary of Proposition~\ref{compact}.
\begin{prop}\label{compactbis}
Let $ s>2 $ be an integer. Then for every $\phi \in H^s_{-1}(\R^2)$
the local solution constructed in Proposition~\ref{compact} can be extended to
a maximal existence interval $[0,T^{\star}[$ such that either
$T^{\star}=\infty$ or
$$
\lim_{t\rightarrow T^{\star}}\|\nabla_{x,y}u\|_{L^{1}_{t}L^{\infty}_{xy}}=\infty.
$$
\end{prop}
\begin{proof}
It suffices to iterate the result of Proposition~\ref{compact} by invoking
(\ref{key}) at each iteration step. 
\end{proof}
It results from Proposition~\ref{compactbis} that the key quantity for the global existence
in $H^{s}_{-1}(\R^2)$, $s>2$ is $\|\nabla_{x,y} u(t,\cdot)\|_{L^{\infty}}$.
\section{A priori estimates using conservation laws} 
In this section we control the growth of some quantities
directly related to the conservation laws of KP I.  Recall that the solutions
obtained in Proposition \ref{compact} satisfy
\begin{eqnarray}\label{kape1}
u_t + u_{xxx} +u u_x +\partial_x (\psi u)-\pa u_{yy}=0.
\end{eqnarray}
In \cite{ZS}, it is shown that the KP-I equation has a Lax pair
representation. This in turn provides an algebraic procedure generating an
infinite sequence of conservation laws. More precisely, if $u$ is a formal
solution of the KP-I equation then
$$
\frac{d}{dt}\Big[\int \chi_n\Big]=0,
$$
where $\chi_1=u$, $\chi_2=u+i\partial_{x}^{-1}\partial_y u$ and for $n\geq 3$,
$$
\chi_n=\Big(\sum_{k=1}^{n-2}\chi_k\, \chi_{n+1-k}\Big)+\partial_x\chi_{n-1}+
i\partial_{x}^{-1}\partial_{y}\chi_{n-1}\,\, .
$$
For $n=3$, we find the conservation of the $L^2$ norm, $n=5$ corresponds to
the energy functional giving the Hamiltonian structure of the KP-I
equation. As we noticed in \cite{MST}, there is a serious analytical
obstruction to give sense of $\chi_9$ as far as $\R^2$ is considered as a
spatial domain. \\

Inspired by the above discussion, we define the following functionals
$$
M(u) =  \int_{\R^2}u^{2}\,, \quad
E(u) =  \frac{1}{2} \intr_{\R^2}u_x^2+\frac{1}{2} \intr_{\R^2} (\pa u_y)^2-
\frac{1}{6}
\intr_{\R^2} u^3\, 
$$
and
\begin{eqnarray*}
F^{\psi}(u) & = &  \frac{3}{2} \intr_{\R^2} u_{xx}^2 + 5 \intr_{\R^2} u_y^2
+ \frac{5}{6}\intr_{\R^2} (\p^{-2}_x u_{yy} )^2
 -\frac{5}{6} \intr_{\R^2} u^2 \p^{-2}_x u_{yy}
\\
& &
-\frac{5}{6} \intr_{\R^2} u\, (\pa u_y)^2
 +\frac{5}{4} \intr_{\R^2} u^2\, u_{xx} + \frac{5}{24} \intr_{\R^2} u^4 
\\
& &
-\frac{5}{3} \intr_{\R^2}\psi\, u\, \p^{-2}_x\, u_{yy}-\frac{5}{6}
\intr_{\R^2}\psi\, 
(\p^{-1}_x u_y )^2\,\, .
\end{eqnarray*}
Recall that the functionals $M$ and $E$ corresponds to the momentum and energy
conservations respectively while the functional $F^{\psi}$ is motivated by the
higher order conservation laws for the KP-I equation associated to $\chi_7$.
Let us notice however that the functional $ F^{\psi}(\cdot) $ contains two
supplementary terms involving $ \psi $ with respect to the
corresponding conservation law of the KP-I equation.
\\

We denote by $H^{\infty}_{-1}(\R^2)$  the intersection of all
$H^{s}_{-1}(\R^2)$. The next proposition gives bounds on the quantities $M$,
$E$ and $F^{\psi}$ for data in spaces where the
local well-posedness of the previous section holds.
\begin{prop} \label{nico1}
For every $R>0$ there exists $C>0$ such that if $ u\in L^{\infty}([0,T];H^{\infty}_{-1}(\R^2))$ is a solution to 
(\ref{integrateLineKP}) corresponding to an initial data $ \phi\in Z\cap
H^{\infty}_{-1}(\R^2)$, $\|\phi\|_{Z}\leq R$ then $E(u(t))$ and $F^{\psi}(u(t))$ are well-defined and $
\forall t\in [0,T ]$,
 \begin{eqnarray}
M(u(t)) & \leq & C\, \exp (C\, t) M(\phi) \label{estL2}\\
  |E(u(t))|&\leq & C\,  \exp( C\, t) |E(\phi)|+g_1(t)  \label{estM}\\
 |F^{\psi}(u(t))|&\leq & C\, \exp( C\, t) |F(\phi)|+g_2(t) \label{estF}
\end{eqnarray}
where $g_1,\, g_2 \, : \R_+\to \R_+ $ are continuous bijections depending only on $R$.
\end{prop}
\begin{proof}[Proof of Proposition~\ref{nico1}.] 
Before entering into the proof of Proposition~\ref{nico1}, we state an
anisotropic Sobolev inequality which will be used in the proof.
\begin{lemma}\label{sobolev}
For $2\leq p\leq 6$ there exists $C>0$ such that for every $u\in
H^{\infty}_{-1}(\R^2)$,
\begin{equation}\label{malak}
\|u\|_{L^p(\R^2)}\leq C\|u\|_{L^2(\R^2)}^{\frac{6-p}{2p}}\,\,
\|u_x\|_{L^2(\R^2)}^{\frac{p-2}{p}}\,\,
\|\partial_x^{-1}u_y\|_{L^2(\R^2)}^{\frac{p-2}{2p}}\,\, .
\end{equation}
\end{lemma}
\begin{proof}
We refer to \cite{BIN} for a proof of (\ref{malak}) and for a systematic study
of anisotropic Sobolev embeddings. 
For a sake of completeness, here we reproduce the proof of (\ref{malak}) given
in \cite{Tom}. Inequality (\ref{malak}) clearly holds for $p=2$. By convexity,
it suffices thus to prove it for $p=6$. 
Following the Gagliardo-Nirenberg proof of the Sobolev embedding, we write
$$
u^{2}(x,y)=2\int_{-\infty}^{x}u_x(z,y)u(z,y)dz
$$
and therefore using the Cauchy-Schwarz inequality in the $z$ integration, we
obtain that for a fixed $y$,
$$
\Big[\sup_{x\in\R}|u(x,y)|\Big]^{4}\leq 4
\Big(\int_{-\infty}^{\infty}u_x^{2}(z,y)dz\Big)\Big(\int_{-\infty}^{\infty}u^{2}(z,y)dz\Big)\, .
$$
Therefore by writing $u^6=u^4 u^2$, we get
\begin{equation}\label{malak-bis}
\int_{\R^2}u^6(x,y)dxdy\leq
4\int_{-\infty}^{\infty}
\Big(
\int_{-\infty}^{\infty}
u_x^{2}(z,y)dz\Big)\Big(
\int_{-\infty}^{\infty}
u^{2}(z,y)dz\Big)^2dy\, .
\end{equation}
Next, using Fubini theorem and an integration by parts, we obtain
\begin{eqnarray*}
\int_{-\infty}^{\infty}u^{2}(z,y)dz & = &
2\int_{-\infty}^{\infty}\int_{-\infty}^{y}u(z,w)u_y(z,w)dwdz
\\
& = & 
2\int_{-\infty}^{y}\int_{-\infty}^{\infty}u(z,w)u_y(z,w)dzdw
\\
& = &
-2\int_{-\infty}^{y}\int_{-\infty}^{\infty}u_{x}(z,w)\,\,\partial_x^{-1}u_y(z,w)dzdw\, .
\end{eqnarray*}
An application of the Cauchy-Schwarz inequality now gives
$$
\sup_{y\in\R}\Big(\int_{-\infty}^{\infty}u^{2}(z,y)dz\Big)^{2}
\leq
4\|u_x\|_{L^2_{xy}}^{2}\|\partial_{x}^{-1}u_y\|_{L^2_{xy}}^{2}\,\, .
$$
Coming back to (\ref{malak-bis}) yields
$$
\int_{\R^2}u^6(x,y)dxdy\leq 16
\|u_x\|_{L^2_{xy}}^{4}\|\partial_{x}^{-1}u_y\|_{L^2_{xy}}^{2}
$$
which is (\ref{malak}) for $p=6$.
This completes the proof of Lemma~\ref{sobolev}.
\end{proof}
Let us return to the proof of Proposition~\ref{nico1}.
Note that since $u$ satisfies   (\ref{integrateLineKP}), one has
$$
u_t\in C([0,T]; H^{s}(\R^2)),\quad \forall\, s\in\N\, .
$$
First, taking the $L^2$-scalar product of \re{integrateLineKP} with $u$, one obtains
$$
\frac{1}{2} \frac{d}{dt} \intR u^2 = \intR \psi u u_x
=-\frac{1}{2} \intR \psi_x u^2
$$
and thus
\begin{equation}\label{l2}
\| u(t)\|_{L^2} \lesssim \exp(Ct\|\psi_x \|_{L^\infty})\, \|\phi\|_{L^2}\,.
\end{equation}
We will use, following \cite{mol}, an exterior regularization of
\re{integrateLineKP} by a sequence of smooth functions $\vp$ that
cut the low frequencies. More precisely,  let $\vp$ be defined via its Fourier
transform as 
\beq \label{b2}
\hat{\vp}(\xi,\eta):= \left\{
\begin{array}{ll}
1 & \hbox{ if }
\ep < |\xi|<\frac{1}{\ep}\hbox{ and }
 \ep < |\eta |<\frac{1}{\ep}, \\
0 & \hbox{otherwise. }
\end{array}
\right. 
\eneq
Note that thanks to the Lebesgue dominated convergence theorem, if $u\in
H^{s}_{-1}(\R^2)$, $s\in\R$ then $\vp\ast u$ converges to $u$ in $H^{s}_{-1}(\R^2)$.
Thanks to the Sobolev embedding similar statements hold for $u\in L^p(\R^2)$,
$2\leq p\leq\infty$.
\\

Taking the convolution of \re{integrateLineKP} with $ \vp $, one gets
\beq 
\vp\ast u_t + \vp \ast u_{xxx} +\vp \ast \partial_x (\psi u
+u^2/2) -\vp \ast \pa u_{yy}=0 \lab{b5} 
\eneq
Setting 
$$ 
u^\ep=\vp\ast u, 
$$ 
multiplying \re{b5} by
$$
-\vp \ast u_{xx} +\p_{x}^{-2} (\vp\ast u_{yy})-\frac{1}{2} (\vp\ast u^2)
$$
and integrating in $ \R^2 $, one obtains that
\begin{multline}\label{energia}
\frac{1}{2}\frac{d}{dt} \Bigl[\int_{\R^2} (\tu_x)^2  +  
 \int_{\R^2} (\partial_x^{-1}\tu_y )^2 -\int_{\R^2} \frac{(\tu)^3}{3}
\Bigr]=
\\
\frac{1}{2}
\int_{\R^2} \Bigr[ \vp \ast u^2-( \tu)^2\Bigl]
 \, \tu_t 
+\intR (\vp_x \ast(\psi u)) u^\ep_{xx}
\\
-\intR (\vp_x\ast (\psi u)) \partial_x^{-2}
u^\ep_{yy} 
+\frac{1}{2} \intR (\vp_x \ast(\psi u)) (\vp \ast u^2) 
\\
:= I+II+III+IV \, .
\end{multline}
Our aim is to passe to the limit $\ep\rightarrow 0$.
Let us first estimate $I$. This argument is very typical for the present
analysis and a similar situation will appear frequently in the rest of the
proof of Proposition~\ref{nico1}. Using equation (\ref{b5}) and the Cauchy-Schwartz
inequality, we can write
\begin{eqnarray*}
|I|& \lesssim &
\|\tu_t \|_{L^2(\R^2)}\| \vp \ast u^2-( \tu)^2\|_{L^2(\R^2)}
\\
& \lesssim &
\Big(\|u\|_{H^3_{-1}(\R^2)}+\|u\|_{H^3_{-1}(\R^2)}^{2}\Big)\| \vp \ast u^2-( \tu)^2\|_{L^2(\R^2)} \,.
\end{eqnarray*}
Next, we write using the triangle inequality and the Sobolev inequality,
\begin{eqnarray*}
\| \vp \ast u^2-( \tu)^2\|_{L^2(\R^2)} 
&\leq & 
\| \vp \ast u^2-u^2\|_{L^2(\R^2)}+\| u^2-( \tu)^2\|_{L^2(\R^2)}
\\
& \leq &
\| \vp \ast
u^2-u^2\|_{L^2(\R^2)}+\|u-u^{\varepsilon}\|_{L^2}(\|u\|_{L^{\infty}}+\|u^{\varepsilon}\|_{L^{\infty}})
\\
&\lesssim &
\| \vp \ast u^2-u^2\|_{L^2(\R^2)}+\|u\|_{H^{2}(\R^2)}\|u-u^{\varepsilon}\|_{L^2}\,.
\end{eqnarray*}
Since $u\in H^{\infty}_{-1}(\R^2)$, we can apply the Lebesgue dominated
convergence theorem to conclude that
$$
\lim_{\ep\rightarrow 0}\| \vp \ast u^2-( \tu)^2\|_{L^2(\R^2)} =0\,.
$$
Therefore the term $I$ tends to zero as $\ep$ tends to zero. 
\\

Next, we write
\begin{eqnarray*}
II & =& \intR (\psi u^\ep)_x u^\ep_{xx} + \intR [\vp\ast (\psi
u)_x -(\psi u^\ep)_x] u^\ep_{xx}
\\
& = & \intR \psi_x u^\ep u^\ep_{xx} +\intR \psi u^\ep_x u^\ep_{xx} +\intR [\vp\ast (\psi u)_x -(\psi u^\ep)_x] u^\ep_{xx}
\\
& =& -\frac{3}{2} \intR \psi_x (u^\ep_x)^2 +\frac{1}{2} \int \psi_{xxx} (u^\ep)^2 +\intR [\vp\ast (\psi u)_x -(\psi u^\ep)_x]
u^\ep_{xx}
\end{eqnarray*}
and
\begin{eqnarray*}
III & =& -\intR (\psi u^\ep)_x \partial^{-2}_x u^\ep_{yy}  -\intR
[\vp\ast (\psi u)_x -(\psi u^\ep)_x] \partial_x^{-2} u^\ep_{yy}
\\
 & = & -\intR  \psi u^\ep_y  \partial^{-1}_x u^\ep_{y}-\intR \psi_y u^\ep \pa u^\ep_y
 -\intR [\vp\ast (\psi u)_y -(\psi u^\ep)_y]
\pa u^\ep_{y}
\\
& =& \frac{1}{2} \intR \psi_x (\pa u^\ep_y)^2 -\intR \psi_y u^\ep \pa u^\ep_y-\intR [\vp\ast
(\psi u)_y -(\psi u^\ep)_y] \pa u^\ep_{y}
\end{eqnarray*}
and
\begin{eqnarray*}
IV&= & \frac{1}{2} \intR  (\psi u^\ep)_x  (u^\ep)^2+\frac{1}{2}
\intR [(\vp_x\ast (\psi u))(\vp \ast u^2)-(\psi u^\ep )_x (u^\ep)^2]
\\
 & =& \frac{1}{3} \intR \psi_x (u^\ep)^3 +\frac{1}{2}
\intR [(\vp_x\ast (\psi u))(\vp \ast u^2)-(\psi u^\ep )_x (u^\ep)^2].
\end{eqnarray*}
Similarly to the analysis for $I$, thanks to the Lebesgue theorem, all
commutator type terms 
involved in $II$, $III$, $IV$
tend to $0$ as $\ep\rightarrow 0$.
Using  Lemma~\ref{sobolev}, we get the bound
$$
\frac{1}{6}\intR |u|^3 \leq C  \|u\|_{L^2}^{3/2} \|u_x\|_{L^2} \|\pa u_y\|_{L^2}^{1/2}
\leq  \frac{1}{4}\|u_x\|_{L^2}^2+
\frac{1}{4}\|\pa u_y\|_{L^2}^2 + C\|u\|_{L^2}^6 \,.
$$
We therefore obtain that
$$
|E(u)|\geq \frac{1}{4}\|u_x\|_{L^2}^2+
\frac{1}{4}\|\pa u_y\|_{L^2}^2 -C\|u\|_{L^2}^6\,.
$$
Integrating (\ref{energia}) on $(0,t)$ for $t\in (0,T]$ gives
\begin{eqnarray*}
|E(u^{\ep}(t))-E(u^{\ep}(0))| & \lesssim & \int_{0}^{t}\Big(
(\|\psi_x \|_{L^\infty}+\|\psi_y\|_{L^\infty})
(|E(u^{\ep}(\tau))|+\|u^{\ep}(\tau)\|^6_{L^2} )
\\
& &
+(\|\psi_{xxx}\|_{L^\infty}+\|\psi_y\|_{L^\infty}) |u^{\ep}(\tau)|_{L^2}^2\,
\Big) d\tau+\int_{0}^{t}|A^{\ep}(\tau)|d\tau,
\end{eqnarray*}
where $\lim_{\ep\rightarrow 0}A_{\ep}(\tau)=0$ for every $\tau\in [0,t]$
($A_{\ep}(\tau)$ corresponds to $I$ and the commutator terms involved in $II$,
$III$, $IV$).
Passing to the limit $ \ep\rightarrow 0 $, we infer that 
\begin{eqnarray*}
|E(u(t))|-|E(\phi)| & \lesssim & 
\int_0^t\Big( (\|\psi_x\|_{L^\infty}+\|\psi_y\|_{L^\infty})
(|E(u(\tau))|+ \|u(\tau)\|^6_{L^2})
\\ 
& & + (\|\psi_{xxx}\|_{L^\infty}+\|\psi_y\|_{L^\infty}) |u(\tau)|_{L^2}^2\,
\Big) d\tau\, .
\end{eqnarray*}
By the Gronwall lemma and (\ref{l2}), it follows that
\begin{equation*}
|E(u(t))| \lesssim \exp( C \,  t) |E(\phi)|+g(t) .
\end{equation*}
where $ g(t) $ is an increasing bijection of $ \R_+ $  which
depends only on $\|\phi\|_{L^2}$.
\\

We now turn to the bound on $ F^\psi(u(t)) $. It is worth noticing that
$$
F^\psi(u)=F(u)-\frac{5}{3} \intr_{\R^2}\psi   u \, \p^{-2}_x u_{yy} - \frac{5}{6} \intr_{\R^2}\psi   (\pa u_y)^2 \quad ,
$$
where $F$ is the corresponding conservation law of the KP-I equation (see \cite{MST}).
The introduction of two additional terms in $F^{\psi}$ is the main new idea in
this paper. Indeed, if we multiply (\ref{b5}) by the multiplier used in
\cite{MST}, we obtain terms which can not be treated as remainders (see $A_1$
and $A_2$ below). The term
$$
-\frac{5}{3} \intr_{\R^2}\psi   u \, \p^{-2}_x u_{yy}
$$
is introduced in the definition of $F^{\psi}$ in order to cancel such ``bad''
remainders. The second additional term in the definition of $F^{\psi}$ is
needed for the proof of the continuous dependence of the flow map on the space $Z$.
\\

As in \cite{MST}, a difficulty in the sequel comes from the fact
that the variational derivative $(F^\psi)'(v)$ contains a term
$c\partial_{x}^{-2}\partial_{yy}(v^2+2\psi v)$. Recall that the formal
derivation of the conservation laws consists in multiplying the KP-I equation
with $F'(u)$ where $u$ is a solution. This procedure meets a difficulty since $\partial_{x}^{-2}$ acts only
on functions with zero $x$ mean value which is not a priori the case of
$v^2+2\psi v $. We overcome this difficulty by introducing the functional $F^{\psi,\ep} $ defined by
\begin{multline*}
F^{\psi,\ep}(u(t)) :=  F^{\psi}(\tu(t)) + \frac{5}{6}\Bigl[  \intr_{\R^2} (\tu(t))^2 \p^{-2}_x \tu_{yy}(t)
- \intr_{\R^2} (\vp \ast u^2(t)) \partial_{x}^{-2}\tu_{yy}(t) \Bigr] 
\\
+\frac{5}{3}\Bigl[\intr_{\R^2}\psi \tu(t)\, \p^{-2}_x \tu_{yy}(t)- \intr_{\R^2} \vp\ast (\psi u (t)) \partial_{x}^{-2}\tu_{yy}(t)
\Bigr] \, .
\end{multline*}
(recall that $ \tu(t):=\vp\ast u(t) $).
We are now in position to state the following lemma.
\begin{lemma}\label{ihp}
Under the assumptions of Proposition~\ref{nico1},
 \begin{eqnarray}
  \frac{d}{ dt}
F^{\psi,\ep}(u(t))  & = & -\frac{5}{3}\intR \psi_x (\partial^{-2}_x
\tu_{yy}) \tu_{xx}
 +\frac{5}{6} \intR \psi (\partial_x^{-2} \tu_{yy})  ((\tu)^2+2 \psi \tu)_x
 \nonumber \\
  & & +\frac{5}{3} \intR \psi \tu (\partial_x^{-2} \tu_{yy})
  +\frac{5}{3} \intR \psi_y(\pa \tu_y) (\partial_x^{-2} \tu_{yy})
\nonumber \\
  & & +G_{\psi}(u^{\ep})+\intR (\partial_x^{-2} \tu_{yy}) \Lambda_\ep(u)+\intR
  \widetilde{\Lambda}_\ep(u)
\end{eqnarray}
 where $ G_{\psi} $ is a continuous functional on  
$$ 
X=\{v\in {\mathcal S}'(\R^2)\,:\,
\|v\|_{L^{2}}+\|\partial_x^{-1} v_y\|_{L^2} +\|v_{xx}\|_{L^2}+\|v_y\|_{L^2} <\infty \, \} 
$$
Moreover
\begin{equation}\label{conv}
\sup_{t\in [0,T]}\intr_{\R^2} |\Lambda_\ep(u(t))|^2
\tendsto{\ep\to 0} 0  \quad{\rm and } \quad \sup_{t\in
[0,T]}\intr_{\R^2} |\widetilde{\Lambda}_\ep(u(t))| \tendsto{\ep\to 0} 0 .
\end{equation}
\end{lemma}
\begin{proof}[Proof of Lemma \ref{ihp}.]
After a direct computation, using (\ref{b5}) and the definition of $
F^{\psi,\ep} $, we obtain the identity
\begin{eqnarray}\label{EqF}
\frac{d}{dt} F^{\psi,\ep}(u) & = & -\frac{5}{3} \, \frac{d}{dt}\intR (\vp \ast (\psi u))(\p^{-2}_x
\tu_{yy})-\frac{5}{6} \frac{d}{dt}\intR \psi (\pa \tu_y)^2
\\
\nonumber
& &
+\intR (A+B) (\tu_{xxx}+ \vp\ast (u u_x) -\pa \tu_{yy})
\\
\nonumber
& &
+\frac{5}{3}\intR \tu (\pa \tu_y) (\tu_{xxy}+\vp\ast (u u_y)
- \partial^{-2}_x \tu_{yyy})
\\
\nonumber
& &
+\frac{5}{3}\intR[u^{\ep}u^{\ep}_{t}-\vp\ast(uu_t)]\p^{-2}_x \u_{yy}
\\
\nonumber
& &
+\intR A  (\vp \ast \partial_x (\psi u))+\intR B (\vp \ast \partial_x (\psi
u)) 
\\
\nonumber
& &
+\frac{5}{3}\intR \tu (\pa \tu_y )(\vp \ast (\psi u)_y)  
\\
\nonumber
& &
:=I+II+III+IV+V+VI+VII+VIII
\end{eqnarray}
where
$$
A=-\frac{5}{3} \, \p^{-4}_x \tu_{4y} +\frac{5}{6}\, \p^{-2}_x \varphi_{\ep,yy}\ast u^2 +\frac{5}{3} \,\u \p^{-2}_x \u_{yy}
$$
and
\begin{eqnarray*}
B= \frac{5}{6}(\pa \tu_y)^2 -\frac{5}{6}(\tu)^3 - 3\tu_{4x} +
10\tu_{yy} -\frac{5}{2}\tu\tu_{xx}
-\frac{5}{4}(\tu)^2_{xx}  \quad .
\end{eqnarray*}
Next one can check that (see \cite[Lemma~1]{MST} for a similar computation)
\begin{equation*}
III+IV+V= \intR (\partial_x^{-2} \tu_{yy}) 
\Lambda^{1}_{\ep}(u)
+\intR\widetilde{\Lambda}_\ep^{1}(u)
\end{equation*}
where 
$$
\Lambda^{1}_{\ep}(u)
=  
-\frac{5}{3}[\vp \ast (u u_t)-\tu \tu_t]
+\frac{5}{3}[\vp \ast (u u_x)-\tu \tu_x]\tu,
$$
and
\begin{eqnarray*}
\widetilde{\Lambda}_\ep^{1}(u)
& =  &
\Bigl[ \vp \ast (u u_x) -\tu \tu_x \Bigr]\times
\\
& &
\times\Bigl( \frac{5}{6}(\pa \tu_y)^2 
-\frac{5}{6}(\tu)^3 
- 3\tu_{4x} +
\frac{25}{3}\tu_{yy}
-\frac{5}{2}\tu\tu_{xx} 
-\frac{5}{4}(\tu)^2_{xx} \Bigr)
\\
& &
+ \frac{5}{3}\Bigl[ \vp \ast (uu_y) - \tu \tu_y \Bigr] \tu \pa \tu_y
\end{eqnarray*}
Using the Lebesgue dominated convergence theorem, we obtain that
$$
\intr_{\R^2} |\Lambda^{1}_{\ep}(u) |^2 \tendsto{\ep\to 0} 0, \quad
\intr_{\R^2} |\widetilde{\Lambda}_\ep^{1}(u) |\tendsto{\ep\to 0} 0, 
\quad t\in [0,T].
$$
Let us next compute the five other terms in the right hand-side of (\ref{EqF}) one by one.
\begin{eqnarray*}
VI= \intR A (\vp \ast \partial_x (\psi u))  & =& \frac{5}{3} \intR (\vp\ast (\psi u)) \partial^{-3}_x \tu_{4y}
\\
& &
-\frac{5}{6} \intR (\vp \ast (\psi u))( \pa \vp_{yy}\ast u^2)  
\\
 &  & + \frac{5}{3} \intR (\vp_x \ast (\psi u)) \tu (\partial_x^{-2} \tu_{yy}) 
\\
& := & A_1+A_2+A_3 \quad .
\end{eqnarray*}
Next
\begin{eqnarray*}
I=-\frac{5}{3} \, \frac{d}{dt}\intR (\vp \ast (\psi u)) 
(\partial^{-2}_x \tu_{yy}) 
& =& -\frac{5}{3} \intR\Big(\partial^{-2}_x\partial_y^2 \vp\ast (\psi u)\Big) \tu_t
\\
& &
-\frac{5}{3} \intR \psi (\partial^{-2}_x \tu_{yy}) \tu_t  
\\
&  &-\frac{5}{3} \intR \psi_t \tu (\partial^{-2}_x \tu_{yy})
\\
& &
-\frac{5}{3} \intR [(\vp\ast (\psi u)_t-(\psi \tu)_t] \partial^{-2}_x \tu_{yy} \nonumber \\
 &  := & C_1+C_2+C_3+C_4  
\end{eqnarray*}
with
\begin{eqnarray*}
C_1
& = & \frac{5}{3} \intR (\partial^{-2}_x \partial_y^2 \vp\ast (\psi u)) \tu_{3x}
-\frac{5}{3} \intR \Big(\partial^{-2}_x \partial_y^2 \vp\ast (\psi u)\Big) \pa \tu_{yy} \\
&  & +\frac{5}{3} \intR (\partial^{-2}_x \partial_y^2 \vp\ast (\psi u))(\vp\ast( u u_x))
+\frac{5}{3} \intR (\partial^{-2}_x \partial_y^2 \vp\ast (\psi u))( \vp\ast \partial_x ( \psi u )) \\
&  =  &\frac{5}{3} \intR  (\partial_y^2 \vp\ast (\psi u)) \tu_{x}-A_1-A_2+0 \quad .
\end{eqnarray*}
As mentioned before, here is the crucial cancellation, thanks to the first
additional term in $F^{\psi}$. Next
\begin{eqnarray*}
\frac{5}{3} \intR  (\partial_y^2 \vp\ast (\psi u)) \tu_{x} & =&\frac{5}{3} \intR \partial_y^2 
(\psi \tu) \tu_{x}+\frac{5}{3} \intR  \Bigl[\partial_y^2 \vp\ast (\psi u) - \partial_y^2  (\psi \tu) \Bigr] \tu_{x} \\
 & =& \frac{5}{6}\intR  \psi_x (\tu_{y})^2+\frac{5}{3} \intR (\psi_{yy} u^\ep +\psi_y u_y^\ep) u_x^\ep \\
& & + \frac{5}{3} \intR  \Bigl[\partial_y^2 \vp\ast (\psi u) - \partial_y^2  (\psi \tu) \Bigr] \tu_{x}
\end{eqnarray*}
Therefore,
\begin{eqnarray*}
VI+C_1 & = & \frac{5}{6}\intR  \psi_x (\tu_{y})^2
+ \frac{5}{3} \intR \psi_{yy} u^\ep u_x^\ep +\frac{5}{3}\intR \psi_y u_y^\ep
u_x^\ep 
\\
& & +\frac{5}{3} \intR (\vp_x \ast (\psi u)) \tu (\partial_x^{-2} \tu_{yy})
\\
& & 
+\frac{5}{3} \intR  \Bigl[\partial_y^2 \vp\ast (\psi u) - \partial_y^2  (\psi \tu) \Bigr] \tu_{x}
\end{eqnarray*}
On the other hand,
\begin{eqnarray*}
C_2 & =& \frac{5}{3} \intR \psi (\partial^{-2}_x \tu_{yy}) \tu_{3x} -\frac{5}{3} \intR
\psi (\partial^{-2}_x \tu_{yy}) \partial^{-1}_x \tu_{yy} \\
 & & +
 \frac{5}{6} \intR \psi (\partial^{-2}_x \tu_{yy})\Big( \vp_x  \ast ((\tu)^2+2\psi \tu )\Big) \\
  & :=& C_{21}+C_{22}+C_{23} \quad .
\end{eqnarray*}
First,
$$
C_{22}=\frac{5}{6} \intR \psi_x (\partial^{-2}_x \tu_{yy} )^2 \quad .
$$
Next,
\begin{eqnarray*}
C_ {21} & =& -\frac{5}{3} \intR
\psi (\partial^{-1}_x \tu_{yy}) \tu_{xx} - \frac{5}{3} \intR \psi_x (\partial^{-2}_x \tu_{yy}) \tu_{xx} \\
 & = & \frac{5}{3} \intR
\psi  \tu_{yy} \tu_{x}+\frac{5}{3} \intR \psi_x  (\pa \tu_{yy}) \tu_{x} -\frac{5}{3} \intR \psi_x (\partial^{-2}_x \tu_{yy}) \tu_{xx} \\
& =& \frac{5}{6} \intR \psi_x (\tu_y)^2 -\frac{5}{3} \intR \psi_x  \tu_{yy} \tu \\
&  &  - \frac{5}{3} \intR \psi_{xx}  (\pa\tu_{yy}) \tu- \frac{5}{3} \intR
\psi_x (\partial^{-2}_x \tu_{yy}) \tu_{xx} 
-\frac{5}{3}\intR \psi_y \tu_y \tu_x
\\
 & = & \frac{15}{6} \intR \psi_x (\tu_y)^2 -\frac{5}{6} \intR \psi_{xxx} (\pa \tu_y)^2
- \frac{5}{3} \intR \psi_x (\partial^{-2}_x \tu_{yy}) \tu_{xx}\\
 & & -\frac{5}{3}\intR \psi_y \tu_y \tu_x -\frac{5}{6}\intR \psi_{xyy}
 (u^\ep)^2+\frac{5}{3}\intR \psi_{xxy} (\pa \tu_y) \, \tu \,.
\end{eqnarray*}
Since $ \psi_t=-\psi_x $,
$$
C_3=\frac{5}{3} \intR \psi_x \tu (\partial_x^{-2} \tu_{yy}) \quad .
$$
Now
\begin{eqnarray*}
-\frac{5}{6} \frac{d}{dt} \intR \psi (\pa \tu_y)^2  & = & 
-\frac{5}{6} \intR \psi_t (\pa \tu_y)^2
\\
& &
+\frac{5}{3} \intR \psi (\pa\tu_y) \tu_{xxy}
\\
& &
+\frac{5}{3} \intR \psi (\pa \tu_y) (\vp\ast (uu_y)) 
\\
& &
+\frac{5}{3} \intR \psi (\pa \tu_y) (\vp\ast (\psi u)_y) 
\\
& &
-\frac{5}{3}
 \intR \psi  (\pa\tu_y)( \partial^{-2}_x \tu_{3y})
\\
&:= &
D_1+D_2+D_3+D_4+D_5\quad .
\end{eqnarray*}
Since $ \psi_t=-\psi_x $ and $ X\hookrightarrow L^\infty(\R^2) $, 
by involving two commutators, we infer that
$$ 
D_1+ D_3 + D_4=G_{1}(\tu)+\intR\widetilde{\Lambda}_\ep^{2}(u)
$$
where $G_1$ is a continuous functional on $X$ and
$$
\intr_{\R^2} |\widetilde{\Lambda}_\ep^{2}(u) |\tendsto{\ep\to 0} 0, 
\quad t\in [0,T].
$$
On the other hand,
\begin{eqnarray*}
D_2 & = & -\frac{5}{3} \intR \psi_x (\pa \tu_y) \tu_{xy}
-\frac{5}{3} \intR \psi \tu_y \tu_{xy} \\
 & = &\frac{5}{3} \intR \psi_{xx} (\pa\tu_y) \tu_y +\frac{5}{3}
 \intR \psi_x (\tu_y)^2+\frac{5}{6} \intR \psi_x (\tu_y)^2 \\
 & =& \frac{15}{6} \intR \psi_x (\tu_y)^2 -\frac{5}{6} \intR
 \psi_{3x} (\pa \tu_y)^2
\end{eqnarray*}
and
\begin{eqnarray*}
D_5 & = & \frac{5}{3} \intR \psi_y (\pa\tu_y)
(\partial^{-2}_x \tu_{yy})
+\frac{5}{3} \intR \psi (\pa \tu_{yy}) (\partial^{-2}_x \tu_{yy}) \\
 & = &\frac{5}{3} \intR \psi_y (\pa\tu_y)
(\partial^{-2}_x \tu_{yy})-\frac{5}{6}
 \intR \psi_x (\partial^{-2}_x \tu_{yy})^2\,.
\end{eqnarray*}
Note that the last term above canceled with $ C_{22} $. Summarizing, we infer that
\begin{eqnarray*}
I+II+VI & = & -\frac{5}{3}\intR \psi_x (\partial^{-2}_x \tu_{yy}) \tu_{xx}
 +\frac{5}{6} \intR \psi (\partial_x^{-2} \tu_{yy})  ((\tu)^2+2 \psi \tu)_x
\\
\nonumber
& &
+\frac{5}{3} \intR (\vp_x \ast (\psi u)) \tu (\partial_x^{-2} \tu_{yy})
 +\frac{5}{3} \intR \psi_x \tu (\partial_x^{-2} \tu_{yy})
\\
\nonumber
& &
+\frac{5}{3} \intR \psi_y(\pa \tu_y)(\partial_x^{-2} \tu_{yy})+G_{2}(\tu)
\\
\nonumber
& &
+\intR (\partial_x^{-2} \tu_{yy}) \Lambda^{2}_{\ep}(u)
+\intR\widetilde{\Lambda}_\ep^{3}(u)
\end{eqnarray*}
where $G_2$ is a continuous functional on $X$ and
$$
\intr_{\R^2} |\Lambda^{2}_{\ep}(u) |^2 \tendsto{\ep\to 0} 0, \quad
\intr_{\R^2} |\widetilde{\Lambda}_\ep^{3}(u) |\tendsto{\ep\to 0} 0, 
\quad t\in [0,T].
$$
Since clearly
$$
VIII=G_{3}(\tu)+\intR\widetilde{\Lambda}_\ep^{4}(u),
$$
where $G_3$ is a continuous functional on $X$ and
$$
\intr_{\R^2} |\widetilde{\Lambda}_\ep^{4}(u) |\tendsto{\ep\to 0} 0, 
\quad t\in [0,T],
$$
it remains to estimate $VII$. We notice that
\begin{eqnarray*}
VII & =&
-\frac{5}{6} \intR (\pa\tu_y)^2 (\vp_x\ast(\psi u)) 
-\frac{5}{6} \intR (\tu)^3 (\vp_x\ast (\psi u))
\\
& &
-3\intR\tu_{4x}(\vp_x \ast (\psi u )) -10\intR \tu_y (\vp\ast (\psi u)_y)
\\
& & 
-\frac{5}{2} \intR \tu \tu_{xx} (\vp_x \ast (\psi u))
-\frac{5}{4} \intR \partial_x^2 (\tu)^2 (\vp_x  \ast(\psi u)) \quad .
\end{eqnarray*}
We now observe that we can write
$$
VII=G_{4}(u^{\ep})+\intR\widetilde{\Lambda}_\ep^{4}(u),
$$
where
$G_4$ is continuous on $X$ and
and
$$
\intr_{\R^2} |\widetilde{\Lambda}_\ep^{4}(u) |\tendsto{\ep\to 0} 0, 
\quad t\in [0,T].
$$
For example,
\begin{eqnarray*}
-3\intR \tu_{4x}(\vp_x \ast (\psi u)) & = & -3\intR \tu_{4x}
(\psi_x \tu +\psi \tu_x) 
\\
& &
-3\intR [\vp\ast (\psi u)_x -(\psi
\tu)_x]\tu_{4x}
\\ & = &
3\intR \tu_{3x}
(2\psi_x \tu_x +\psi_{xx} \tu +\psi \tu_{xx}) 
\\ & &
-3\intR [\vp\ast
(\psi u)_x -(\psi
\tu)_x]\tu_{4x}\\
 &= & -\frac{15}{2} \intR \psi_x |\tu_{xx}|^2
 +\frac{9}{2}\intR\psi_{3x} |\tu_x|^2
 -3\intR \psi_{3x} \tu \tu_{xx} \\
  & & -3\intR [\vp\ast
(\psi u)_x -(\psi \tu)_x]\tu_{4x}
\end{eqnarray*}
All other terms in the representation of $VII$ can be treated similarly.
This achieves the proof of Lemma~\ref{ihp}. 
\end{proof}
Now, since
$$
\Bigl|\frac{5}{6} \intr_{\R^2} [\varphi^{\ep}\ast(u^2+2\psi u)] \p^{-2}_x \tu_{yy} \Bigr|
\le C (\|u \|_{L^4}^4+\|u\|_{L^2}^2 \|\psi\|^2_{L^{\infty}}) + \frac{5}{12} \|\p^{-2}_x \tu_{yy} \|^2_{L^2},
$$
there exists a constant $ C>0 $ such that for $ \ep $ small enough,
\beq
F^{\psi,\ep} (u(t))\geq \frac{5}{24} \|\partial^{-2}_{x} \tu_{yy}(t) \|_{L^2}^{2} -C , \quad
 \forall t\in [0,T] \, .
\label{b16}
\eneq
We thus deduce from Lemma \ref{ihp}  and \re{b16} that
\begin{equation}
\frac{d}{dt} F^{\psi,\ep}(u(t)) \lesssim \|\psi\|_{W^{3,\infty}}
|F^{\psi,\ep}(u(t))| +R_{\psi}(\tu)+\Lambda_\ep(t)
\end{equation}
where $R_{\psi}$ is continuous on $X$ and $ |\Lambda_\ep|_1 \to 0 $ as $ \ep \to 0 $ uniformly for $ t\in [0,T] $. 
Here we used that thanks to Lemma~\ref{sobolev},
$$
\|\pa \tu_y \|_{L^4}\lesssim \|\pa \tu_y \|^{1/4}_{L^2} \|\tu_y\|^{1/2}_{L^2} \|\partial_x^{-2} \tu_{yy} \|_{L^2}^{1/4}
$$
and
$$
\|\tu_x \|_{L^4}\lesssim \| \tu_x \|^{1/4}_{L^2} \|\tu_{xx}\|^{1/2}_{L^2} \| \tu_{y} \|_{L^2}^{1/4}
$$
and thus,
\begin{eqnarray*}
\intR \Bigl| \psi \tu_x \tu (\partial_x^{-2} \tu_{yy}) \Bigr| 
& \lesssim & 
\|\psi\|_{L^\infty} \|\tu\|_{L^4} \|\tu_x\|_{L^4} 
\|\partial_x^{-2} \tu_{yy}\|_{L^2} \\
 & \lesssim & |F^{\psi,\ep}(u)| +R_{\psi}(\tu)
\end{eqnarray*}
and in the same way
\begin{eqnarray*}
\intR \Bigl| \psi  (\partial_x^{-1} \tu_{y}) \tu \tu_y   \Bigr| & \lesssim &
\|\psi\|_{L^\infty} \|\tu\|_{L^4} \|\tu_x\|_{L^4} 
\|\partial_x^{-1} \tu_{y}\|_{L^4}  \|\tu_y\|_{L^4} \\
 & \lesssim & |F^{\psi,\ep}(u)| +R_{\psi}(\tu)
\end{eqnarray*}
Hence,
\begin{equation}
|F^{\psi,\ep}(u(t))|\lesssim \exp(Ct )|F^{\psi,\ep}(\phi)| +
\exp(Ct). \label{est7}
\end{equation}
Letting $ \ep $ tends
to $ 0 $, $
F^{\psi,\ep} ( \phi) \to
 F^\psi (\phi)< \infty $ and thus
$$
\sup_{t\in [0,T], \, \ep>0 } F^{\psi,\ep} (u(t)) \lesssim 1.
\lab{b16b}
$$
>From \re{b16}, one infers
$$
\sup_{t\in [0,T], \, \ep>0} \Bigl| \vp\ast \p^{-2}_x u_{yy}(t)
 \Bigr|_2 \lesssim 1 
$$
 and thus
$$
\p^{-2}_x  u_{yy} \in L^\infty(0,T; L^2(\R^2))  \; ,
$$
by the Lebesgue theorem. It is then easy to check that
$$
F^{\psi,\epsilon}(u(t))
\tendsto{\ep\to 0}
F^\psi(u(t)),\quad t\in [0,T].
$$
and thus \re{estF} follows from \re{est7}. 
This completes the proof of Proposition~\ref{nico1}.
\end{proof} 
It follows from \re{b16} and the Lebesgue theorem that,
for $ u\in C([0,T];H^\infty_{-1}(\R^2)) $, 
$\partial_x^{-2}\tu_{yy} $ tends to $ \partial_x^{-2} u_{yy}  $ in
$   L^\infty([0,T];L^2(\R^2)) $. Therefore we can pass to  the limit in $\ep $ in
Lemma~\ref{ihp} which leads to the next Proposition.
\begin{prop}\label{BorneZ}
Let $ u\in C([0,T];H^\infty_{-1}(\R^2)) $ then
\begin{equation}
\|u(t)\|_Z \lesssim g(t)\quad , \label{borneZ}
\end{equation}
where $ g$ is an increasing continuous one to one mapping of $
\R_+ $ only depending on $ \|u(0)\|_{Z} $. Moreover,
\begin{eqnarray}
F^{\psi}(u(t))-F^{\psi}(u(0))& = & -\frac{5}{3}\int_0^t \intR \psi_x
(\partial^{-2}_x u_{yy}) u_{xx}
\nonumber
\\
& & +\frac{5}{6} \int_0^t \intR \psi (\partial_x^{-2} u_{yy})  ((u)^2+2 \psi u)_x
 \nonumber 
\\
& & +\frac{5}{3} \int_0^t\intR \psi u (\partial_x^{-2} u_{yy})
+\frac{5}{3}\intR(\psi u)_{x}u(\partial_x^{-2} u_{yy})
\nonumber
\\
&&
+\frac{5}{3}\int_0^t \intR \psi_y(\pa u_y)(\partial_x^{-2} u_{yy}) +
\int_0^t G_{\psi}(u)
 \; , 
\label{equF}
\end{eqnarray}
where $G_{\psi}$ is continuous on $X$.
\end{prop}
%
\section{Dispersive estimates}
Let us first extend some linear dispersive estimates established in \cite{Sa} for anti-derivatives in $ x$ of
the free group of KP. A similar argument was used in \cite{MR}. 
Let 
$$
U(t):=\exp(-t(\partial_{x}^{3}-\partial_{x}^{-1}\partial_y^2))
$$ 
be the
unitary group on $H^s(\R^2)$ defining the free KP-I evolution.
Let $D_x$ be the Fourier multiplier with symbol $|\xi|$. Then we have the
following Strichartz inequality for the free KP-I evolution.
\begin{lemma}\label{Strichartz}
Let $T>0$.
Then for every $ 0\leq \varepsilon\le 1/2 $, we have the estimates 
\begin{equation}\label{hom}
\|D_x^{-\frac{\ep \delta(r)}{2}} U(t)\phi  \|_{L^{q}_{T}L^{r}_{xy}}
\lesssim  \|\phi\|_{L^2},
\end{equation}
and
\begin{equation}\label{non-hom}
\left\|\int_{0}^{t}D_x^{-\frac{\ep\delta(r)}{2}}U(t-t')F(t')dt'\right\|_{L^{q}_{T}L^{r}_{xy}}
\lesssim \|F\|_{L^{1}_{T}L^{2}_{xy}},
\end{equation}
provided $r\in [2,\infty]$, and
$$
0\leq \frac{2}{q}=(1-\ep/3)\,\delta(r)<1\, 
$$
with
$$
\delta(r):=1-\frac{2}{r}\,\, .
$$
\end{lemma}
\begin{proof}
Since $U(t)$ and $D_x$ commute, using the
Minkowski inequality, we obtain that (\ref{non-hom}) follows from (\ref{hom}).
Let us now turn to the proof of (\ref{hom}). For $r=2$, (\ref{hom}) follows
from the unitarity of $U(t)$ on $L^2(\R^2)$. It suffices thus to prove it for $r=\infty$
Let us set
$$ 
G(x,y,t):=\int_{\R^2} e^{it(x\xi+y \eta)}\, e^{it(\xi^3-\eta^2/\xi)}\, d\xi
\, d\eta \, .
$$
Then 
$$
U(t)\phi =G(\cdot,\cdot,t)\star \phi\, ,
$$
where $\star$ denotes the convolution with respect to the spatial variables
$x$ and $y$.
Integrating a gaussian integral with respect to $\eta$ (see \cite{Sa}), we get
for $ 0\leq \varepsilon\le 1/2$,
\begin{equation*} 
\|D_x^{-\ep}\, G(x,y,t) \|_{L^\infty_{xy}}\lesssim   |t|^{-1/2}
\sup_{x\in\R}\,\, \Bigl|\int_{0}^{\infty}|\xi|^{\frac{1}{2}-\ep}\, e^{i(t\xi^3+x\xi)}d\xi\Bigr|\, .
\end{equation*}
Next using the Van der Corput lemma as in \cite{KPV}, we infer that
\begin{eqnarray*}
\|D^{-\ep}_x G(x,y,t) \|_{L^\infty_{xy}} \lesssim & |t|^{-1+\ep/3}
\end{eqnarray*}
and thus
\begin{equation}\label{disp-in}
\|D^{-\varepsilon}_{x} U(t)\, \phi \|_{L^\infty_{xy}}\lesssim |t|^{-1+\ep/3}\|\phi \|_{L^1_{xy}}
\end{equation}
which is the key dispersive inequality.
Let us now perform the duality argument which provides (\ref{hom}) for
$r=\infty$ as a consequence of the dispersive inequality (\ref{disp-in}).
We first fix the $q$ corresponding to $r=\infty$, i.e.
$$
\frac{2}{q}=1-\frac{\ep}{3}\, .
$$
Set
$$
A:=D_x^{-\ep/2}U(t)\, .
$$
Our goal is to show that $A$ is bounded from $L^2_{xy}$ to
$L^{q}_{T}L^{\infty}_{xy}$. Denote by $A^{\star}$ the operator formally
adjoint to $A$. Then  
$$
AA^{\star}(f)=\int_{0}^{T}D_x^{-\ep}U(t-t')f(t')dt'\, .
$$
Using the dispersive estimate (\ref{disp-in}) and the Hardy-Littlewood
inequality in time, we infer that
\begin{equation}\label{dual}  
\|AA^{\star}(f)\|_{L^q_{T}L^{\infty}_{xy}} \lesssim \|f\|_{L^{q'}_T L^{1}_{xy}},
\end{equation}
where $q'$ is the conjugate of $q$, i.e. 
$$
\frac{1}{q}+\frac{1}{q'}=1\, .
$$
Recall that $L^{\infty}_{xy}$ can be seen as the dual space of $L^{1}_{xy}$.
It is not true that $L^{q'}_T L^{1}_{xy}$ can be seen as the dual of
$L^q_{T}L^{\infty}_{xy}$ but however, it is easy to see that for $q<\infty$, we can write 
\begin{equation}\label{pisha}
\|Af\|_{L^q_{T}L^{\infty}_{xy}}=
\sup_{\|g\|_{L^{q'}_T L^1_{xy}}\leq 1}|\langle Af\, ,\,  g\rangle|
\end{equation}
where $\langle \cdot,\cdot\rangle$ denotes the $L^2_{T}L^2_{xy}$ inner product.
Let us next write by the Cauchy-Schwarz inequality
\begin{equation}\label{pisha-bis}
|\langle Af\, ,\,  g\rangle|=|(f\, ,\, A^{\star}g)|\leq
\|f\|_{L^2_{xy}}\|A^{\star}g\|_{L^2_{xy}},
\end{equation}
where $(\cdot,\cdot)$ denotes the $L^2_{xy}$ inner product.
Next, using (\ref{dual}), we obtain
$$
\|A^{\star}g\|_{L^2_{xy}}^{2}=(A^{\star}g\, ,\, A^{\star}g)=
\langle AA^{\star}(g)\, ,\, g\rangle
\leq \|AA^{\star}(g)\|_{L^q_{T}L^{\infty}_{xy}}
\|g\|_{L^{q'}_{T}L^{1}_{xy}}
\leq C\|g\|_{L^{q'}_{T}L^{1}_{xy}}^{2}\, .
$$
Coming back to (\ref{pisha}) and (\ref{pisha-bis}) ends the proof of (\ref{hom}) for $r=\infty$.
This completes the proof of Lemma~\ref{Strichartz}. 
\end{proof}
We next state the crucial dispersive estimate.
\begin{prop}\label{carlos}
Let $ v\in C([0,T];H^{\infty}_{-1}(\R^2))$ be a solution of 
\begin{equation}\label{Eqnonhom}
v_t+v_{xxx} -\partial_x^{-1} v_{yy} = F_x \quad .
\end{equation}
Then for every $\ep>0$ there exists $C_{\ep}$ such that 
\begin{equation}\label{estim1}
\| v\|_{L^{1}_T L^\infty_{xy}} \leq C_{\ep}(1+T) 
\|J_x^{1/2+\ep} \,v\|_{L^\infty_T\, L^2_{xy}}
+C_{\ep}\|J_x^{1/2+\ep}F\|_{L^1_{T}L^2_{xy}}\quad .
\end{equation}
\end{prop}
\begin{proof}
We consider a Littlewood-Paley decomposition in the $x$-variable
$$
v= \sum_{\lambda\,-{\rm dyadic}} v_\lambda
$$
where $ v_\lambda := \Delta_{\lambda} v $ and the Fourier multipliers $ \Delta_\lambda $ are defined as follows:
$$
\widehat{\Delta_\lambda v}(t, \xi ,\eta) := \left\{
\begin{array}{ll}
\varphi(\frac{\xi}{\lambda}) \, \hat{v}(t,\xi,\eta), & \lambda=2^k, \, k\ge 1 , \\
 \chi(\xi) \, \hat{v}(t,\xi,\eta) , & \lambda=1,
\end{array}
\right.
$$
where the nonnegative functions $ \chi \in C_0^\infty(\R) $ and $  \varphi \in
C_0^\infty(\R^*)$ are defined as in \cite{KT1}.
For $ \lambda\ge 2 $ fixed we write a natural splitting 
$$ 
[0,T] =\bigcup_{j} I_j 
$$ 
where $I_j=[a_j,b_j]$ are with disjoint interiors and $ |I_j|\leq
\lambda^{-1}$.
Clearly, we can suppose that the number of the intervals $I_j$ is bounded by $C(1+T)\lambda$.   
Using the H\"older inequality in time, we can write
\begin{equation*}
\| v_\lambda \|_{L^{1}_T L^\infty_{xy}} \lesssim  \sum_j \|v_\lambda \|_{L^{1}_{I_j} L^\infty_{xy}} 
\lesssim  \lambda^{-\frac{1}{2}-\frac{\ep}{6}}\,
  \sum_j \|v_\lambda \|_{L^{q_\ep}_{I_j} L^\infty_{xy}}\, ,
\end{equation*}
where $ 1/q_\ep=1/2-\ep/6 $.
Next, we apply the Duhamel formula in each $ I_j $ to obtain on $ I_j $,
$$
v_\lambda(t)=D_x^{-\ep/2} U(t-a_j)D_x^{\ep/2} v_\lambda (a_j) -
\int_{a_j}^tD_x^{-\ep}U(t-t') [\Delta_\lambda D_x^\ep
\partial_x F](t') \, dt'
$$
Using the Strichartz estimates established in Lemma~\ref{Strichartz}, it yields
$$
\|v_\lambda\|_{L^{q_\ep}_{I_j} L^\infty_{xy}} \lesssim
\|D^{\ep/2}_x v_\lambda(a_j)\|_{L^2} + \|\Delta_\lambda
D_x^{1+\ep} F \|_{L^1_{I_j} L^2_{xy}}\, .
$$
Therefore
$$
\|v_\lambda\|_{L^{1}_{I_j} L^\infty_{xy}} \lesssim \lambda^{-1/2+\ep/3} \| v_\lambda(a_j)\|_{L^2} +
 \lambda^{1/2+5\ep/6}\|\Delta_\lambda  F \|_{L^1_{I_j} L^2_{xy}}
$$
and summing over $j$,
\begin{multline}\label{eq5}
\| v_\lambda\|_{L^{1}_{T} L^\infty_{xy}} \lesssim \lambda^{-1/2+\ep/3}\sum_j  \| v_\lambda\|_{L^\infty_T
L^2_{xy}} + \lambda^{1/2+5\ep/6} \|\Delta_\lambda F\|_{L^1_{T} L^2_{xy}}  
\\
\lesssim \lambda^{1/2+\ep/3} (1+T) \| v_\lambda\|_{L^\infty_T
L^2_{xy}} + \lambda^{1/2+5\ep/6} \|\Delta_\lambda
F\|_{L^1_{T} L^2_{xy}} \, .
\end{multline}
Moreover, again by Duhamel formula  and Strichartz estimates
\begin{equation}\label{bassefreq}
\|{\Delta}_1 v \|_{L^{1}_T L^\infty} \lesssim (1+T)\Bigl(\|{\Delta}_1
v(0)\|_{L^2_{xy}}+\|{\Delta}_1 F \|_{L^1_T L^2_{xy}}\Bigr) 
\end{equation}
Hence, by Minkowski and Bernstein inequalities and \re{eq5}-\re{bassefreq},
for any 
$$ 
0<\alpha,\ep<\!\! <  1, 
$$ 
we get
\begin{eqnarray*}
\|v\|_{L^{1}_T L^\infty_{xy}} 
& \lesssim &  \sum_{\lambda}\|v_\lambda\|_{L^{1}_T L^\infty_{xy} }  
\\
& \lesssim & \|{\Delta}_1 v\|_{L^{1}_T L^\infty_{xy}}+
\sum_{\lambda\ge 2} \lambda^{-\alpha}
\lambda^{\alpha}
\|v_\lambda\|_{L^{1}_{T}
L^\infty_{xy}} 
\\
& \lesssim &  \|{\Delta}_1 v\|_{L^{1}_T L^\infty_{xy}} +
\sup_{\lambda\ge 2}\lambda^{\alpha}
\| v_\lambda\|_{L^{1}_{T} L^\infty_{xy}} 
\\
& \lesssim &  (1+T) \|J_x^{\alpha +1/2+\ep/3} \, v\|_{L^\infty_T
 L^2_{xy}}+\|J_x^{\alpha +1/2+5\ep/6}F\|_{L^1_{T}L^2_{xy}}\, .
\end{eqnarray*}
This completes the proof of Proposition~\ref{carlos}.
\end{proof}
\section{A first global existence result}
In this section, we use the dispersive estimates of the previous section and
the a priori bounds of Section~3 to show that the local solutions with
$H^{\infty}_{-1}(\R^2)\cap Z$ data, obtained in Section~2 are in fact {\bf global} in time.
\begin{prop}\label{globexist}
Let $\phi\in H^{\infty}_{-1}(\R^2)\cap Z$. Then there exists a
unique global solution $u$ of (\ref{integrateLineKP}) with
initial data $\phi$, such that 
$$
u\in C(\R_{+};H^{\infty}_{-1}(\R^2))\cap L^{\infty}_{loc}(\R_{+};Z),
\quad
u_{t}\in C(\R_{+};H^{\infty}_{-1}(\R^2))\, .
$$
\end{prop}
\begin{proof}
From the local theory of Proposition~\ref{compact} and Proposition~\ref{compactbis} to prove that
the local solution $u$ can be extended on $ \R_+ $,  it suffices
to show that $ \|\nabla_{x,y} u(t) \|_{L^\infty} $  can not go to
infinity in finite time. From the propagation of regularity obtained in
Proposition~\ref{compact} and Proposition~\ref{BorneZ}, we deduce that $u$ is
bounded in $Z$ as long as it exists. 
Unfortunately $\|\nabla_{x,y}u\|_{L^{\infty}}$ {\bf is not} controlled by
$\|u\|_{Z}$ for an arbitrary function $u$. But thanks to the key
Proposition~\ref{carlos} will be able to control $\|\nabla_{x,y}u\|_{L^{\infty}}$
whenever $u$ is a solution of (\ref{integrateLineKP}). This is the
crucial dispersive effect in our analysis.
\\

According to \cite{KP} and (\cite{KPV2}, Theorem A.12),
 we have the following generalizations of the Leibniz rule.
\begin{lemma} \label{Leibniz}
Let  $ 1<p<\infty $. Then for $ s>0 $,
\begin{equation*}
\|J^s (fg)\|_{L^p(\R^n)}\lesssim \|J^s f\|_{L^{p}(\R^n)}\|g\|_{L^{\infty}(\R^n)}+\| f\|_{L^{\infty}(\R^n)} \|J^s g\|_{L^{p}(\R^n)}
\end{equation*}
and for $ 0<s<1 $,
\begin{equation*}
\|J^s (fg)\|_{L^p(\R)}\lesssim \|J^s f\|_{L^{\infty}(\R)}\|g\|_{L^{p}(\R)}+\| f\|_{L^{\infty}(\R)} \|J^s g\|_{L^{p}(\R)} \quad .
\end{equation*}
\end{lemma}
Proposition~\ref{carlos} yields the following estimates on smooth solutions to \re{kape1}.
\begin{lemma}\label{key-lemma}
For every $ 0<\ep <\!<1 $ there exists $C_{\ep}$ such that if
$ u\in C([0,T];H^{\infty}_{-1}(\R^2))$ is a  solution to (\ref{integrateLineKP}) then
\begin{equation}
\| u_x\|_{L^{1}_T L^\infty_{xy}} \lesssim C_{\ep}\Bigl( 1+T+
\|u\|_{L^1_T L^\infty_{xy}} +T
\big(\|J_x^{\frac{3}{2}+\ep}\psi\|_{L^{\infty}}+\|\psi\|_{L^{\infty}}\big)
\Bigr)
 \| J_x^{\frac{3}{2}+\ep} u\|_{L^\infty_T L^2_{xy}}
\label{estim2}
\end{equation}
and
\begin{eqnarray}
\nonumber
\| u_y\|_{L^{1}_T L^\infty_{xy}}& \lesssim &\,\, C_{\ep}(1+T)\| J_x^{\frac{1}{2}+\ep}u_y\|_{L^\infty_T L^2_{xy}}
\\
\nonumber
&  & +\,
C_{\ep}T\big(\|J_x^{\frac{1}{2}+\ep}\psi_y\|_{L^{\infty}}+\|\psi_y\|_{L^{\infty}}\big)\|J_x^{\frac{1}{2}+\ep}u\|_{L^{\infty}_{T}L^{2}_{xy}}
\\
\label{estim3}
&  & +\,
C_{\ep}T\big(\|J_x^{\frac{1}{2}+\ep}\psi\|_{L^{\infty}}+\|\psi\|_{L^{\infty}}\big)\| J_x^{\frac{1}{2}+\ep}u_y\|_{L^\infty_T L^2_{xy}}
\\
\nonumber
&  & +\,
C_{\ep}\big(\|u\|_{L^1_T L^\infty_{xy}}+\|J_{x}^{1/2+\ep}\,u\|_{L^1_{T}L^{\infty}_{xy}}\big)
\| J_x^{\frac{1}{2}+\ep}u_y\|_{L^\infty_T L^2_{xy}}\,\, .
\end{eqnarray}
\end{lemma}
\begin{proof}
We have that $ u_x $ satisfies \re{Eqnonhom} with 
$$
F=\partial_x (\psi u+u^2/2)\,.
$$
Estimate \re{estim2} follows from \re{estim1} and
Lemma~\ref{Leibniz} applied in the $x$ variable.
More precisely, for a fixed $t,y$, we have that
$$
\|J_x^{\frac{1}{2}+\ep}\partial_x(\psi u)\|_{L^2_x}=
\|J_x^{\frac{1}{2}+\ep}(\psi_x u+\psi u_x )\|_{L^2_x}\lesssim
(\|\psi\|_{L^\infty_x} +\|J_x^{\frac{3}{2}+\ep}\psi\|_{L^{\infty}_x})\|J_x^{\frac{3}{2}+\ep} u\|_{L^2_x}
$$
where we use that $\|J_x^{\frac{1}{2}+\ep}\psi\|_{L^{\infty}_x} \lesssim \|\psi\|_{L^\infty_x} +\|J_x^{\frac{3}{2}+\ep}\psi\|_{L^{\infty}_x} $ (see Lemma \ref{luca} below).
Taking the $L^{1}_{T}L^{2}_y$ norm of the last inequality gives
$$
\|J_x^{\frac{1}{2}+\ep}\partial_x(\psi u)\|_{L^1_{T}L^{2}_{xy}}
\lesssim T
\big(\|J_x^{\frac{3}{2}+\ep}\psi\|_{L^{\infty}}+\|\psi\|_{L^{\infty}}\big)
\| J_x^{\frac{3}{2}+\ep} u\|_{L^\infty_T L^2_{xy}}\, .
$$
Further, for a fixed $t,y$, Lemma \ref{Leibniz} yields
$$
\|J_x^{\frac{1}{2}+\ep}\partial_x(u^2/2)\|_{L^2_x}\lesssim 
\|u\|_{L^{\infty}_x}
\|J_x^{\frac{3}{2}+\ep} u\|_{L^2_x}\, .
$$
Taking the $L^{1}_{T}L^{2}_y$ norm gives
$$
\|J_x^{\frac{1}{2}+\ep}\partial_x(u^2/2)\|_{L^1_{T}L^2_{xy}}\lesssim
\|u\|_{L^1_{T}L^{\infty}_{xy}}
\|J_x^{\frac{3}{2}+\ep} u\|_{L^{\infty}_{T}L^2_{xy}}\, .
$$
This proves that \re{estim2} follows from \re{estim1}.
We next estimate $u_y$. Notice that $u_y$ satisfies \re{Eqnonhom}
with
$$
F=\partial_{y}(\psi u+u^2/2)\, .
$$
Using Lemma~\ref{Leibniz}, we obtain
\begin{eqnarray*}
\|J_x^{1/2+\ep}\partial_y(\psi u+u^2/2)\|_{L^1_T L^{2}_{xy}} & = &
\|J_x^{1/2+\ep}(\psi_y u+\psi u_y+uu_y)\|_{L^1_T L^{2}_{xy}}
\\
&\lesssim & 
T\big(\|J_x^{\frac{1}{2}+\ep}\psi_y\|_{L^{\infty}}+\|\psi_y\|_{L^{\infty}}\big)
\|J_x^{\frac{1}{2}+\ep}u\|_{L^{\infty}_{T}L^{2}_{xy}}
\\
& & 
+T\big(\|J_x^{\frac{1}{2}+\ep}\psi\|_{L^{\infty}}+\|\psi\|_{L^{\infty}}\big)
\|J_x^{\frac{1}{2}+\ep}u_y\|_{L^{\infty}_{T}L^{2}_{xy}}
\\
& &
+
\|u\|_{L^1_T L^{\infty}_{xy}}\|J_x^{\frac{1}{2}+\ep}u_y\|_{L^{\infty}_{T}L^2_{xy}}
+
\|J_x^{\frac{1}{2}+\ep} u\|_{L^1_T L^{\infty}_{xy}}\|u_y\|_{L^{\infty}_{T}L^2_{xy}}
\end{eqnarray*}
which implies \re{estim3}.
This completes the proof of Lemma~\ref{key-lemma}.
\end{proof}
In order to control $\|J_x^{\frac{1}{2}+\ep} u\|_{L^1_T L^{\infty}_{xy}}$ we
will need the following lemma.
\begin{lemma}\label{luca}
For every $0<s<1$ there exists $C$ such that for every $u\in L^{\infty}_{xy}$
satisfying $u_x\in L^{\infty}_{xy}$ one has the estimate
$$
\|J^s_x u\|_{L^{\infty}_{xy}}
\leq
C\big(\|u\|_{L^{\infty}_{xy}}+
\|u_x\|_{L^{\infty}_{xy}}\big)\, .
$$
\end{lemma}
\begin{proof}
Consider a Littlewood-Paley decomposition with respect to the $x$ variable
$$
u=\sum_{\lambda-{\rm dyadic}}\Delta_{\lambda}\, u\, .
$$
We have that $J^{s}_{x}\Delta_1$ is bounded on $L^{\infty}_{xy}$. Moreover one
can directly check that there exists $C>0$ such that for every $\lambda\geq 2$, 
$$
\lambda\|\Delta_{\lambda}u\|_{L^{\infty}_{xy}}\leq
C\|u_x\|_{L^{\infty}_{xy}}\, .
$$
Therefore using the triangle inequality and the Bernstein inequality, we obtain
that
\begin{eqnarray*}
\|J^s_x u\|_{L^{\infty}_{xy}}
& \leq & \| J^s_x \Delta_1 u\|_{L^{\infty}_{xy}}+\sum_{\lambda\geq 2}
\| J^s_x \Delta_{\lambda} u\|_{L^{\infty}_{xy}}
\\
& \lesssim &
\| u\|_{L^{\infty}_{xy}}+\sum_{\lambda\geq 2}\, \lambda^{s}\| \Delta_{\lambda}
u\|_{L^{\infty}_{xy}}
\\
& \lesssim & 
\| u\|_{L^{\infty}_{xy}}+\sum_{\lambda\geq 2}\,
\lambda^{s-1}\|u_x\|_{L^{\infty}_{xy}}
\\
& \lesssim &
\|u\|_{L^{\infty}_{xy}}+ \|u_x\|_{L^{\infty}_{xy}}\,\, .
\end{eqnarray*}
This completes the proof of Lemma~\ref{luca}.
\end{proof}
Let us now return to the proof of  Proposition~\ref{globexist}.
Noticing that $\|u\|_{L^{\infty}}$ is controlled by $\|u\|_{Z}$,  
using (\ref{estim2}) with $\ep=1/2$ and Proposition~\ref{BorneZ}, we obtain
the existence of an increasing continuous function\footnote{The important
point is that this function does not go to infinity in finite time.} $g_1$ from $\R_{+}$ to $\R_{+}$ such that
\begin{equation}\label{u_x}
\|u_x\|_{L^{1}_{T}L^{\infty}_{xy}}\leq g_{1}(T)
\end{equation}
as long as $u(t)$ exists.
Notice that, thanks to Lemma~\ref{luca}, estimate (\ref{u_x}) also provides a
bound for 
$$
\|J^s_x u\|_{L^{1}_{T}L^{\infty}_{xy}},\quad 0<s<1\,\, .
$$
Next, using (\ref{kando1}) and (\ref{kando2}) and the Gronwall lemma, we
obtain that for every integer $s$
\begin{equation}\label{x}
\|\partial_x^{s}u(t,\cdot)\|_{L^2}\leq \|\partial_x^{s}\phi\|_{L^2}
\exp\Big(\|u_x\|_{L^{1}_{t}L^{\infty}_{x,y}}+t\|\psi\|_{W^{s,\infty}}\Big)\, .
\end{equation}
Therefore using Proposition~\ref{BorneZ} to bound $\|u(t,\cdot)\|_{L^2}$ we
obtain that there exists an increasing continuous function $ g_2 $ from $ \R_+ $ to
$ \R_+ $ such that 
\begin{equation}\label{4}
\|J_x^{4} u (t,\cdot)\|_{ L^2_{xy}}\le g_2(t) 
\end{equation}
as long as $u(t)$ exists.
Noticing that
\begin{eqnarray*}
\|J_x^{1/2+\ep}u_y\|_{L^2}^{2}
& = & \|u_y\|_{L^2_{xy}}+
\||\xi|^{1/2+\ep}\, \eta\, \hat{u}(\xi,\eta)\|_{L^2_{\xi\eta}}^{2}
\\
& = & \|u_y\|_{L^2_{xy}}+\||\xi|^{3/2+\ep}\, |\xi|^{-1}\eta\, \hat{u}(\xi,\eta)\|_{L^2_{\xi\eta}}^{2}
\\
& \lesssim & \|u_y\|_{L^2_{xy}}+\|J_x^{4} u\|_{L^2_{xy}}\|\partial_x^{-2}u_{yy}\|_{L^2_{xy}}
\end{eqnarray*}
we deduce from (\ref{4}), (\ref{estim3}), (\ref{u_x}), Lemma~\ref{luca} and Proposition~\ref{BorneZ}
that there exists an increasing continuous function $g_3$ from $ \R+ $ to $ \R_+ $ such that
\begin{equation*}
\|u_y\|_{L^{1}_{T}L^\infty_{xy}} \lesssim g_3(t) 
\end{equation*}
Therefore $ \|\nabla_{x,y} u\|_{L^1_{T}L^\infty_{xy}} $  can not go to
infinity in finite time $T$.
Thanks to the well-posedness Proposition~\ref{compact} 
we deduce that $u$ can be extended on the whole real axis
which completes the proof of Proposition~\ref{globexist}.
\end{proof}
\section{Well-posedness in $ H^{2,0}(\R^2) $}
The next two sections are devoted to the global well-posedness of
(\ref{LineKP}) in $Z$, i.e. we remove the condition $\phi\in H^{\infty}_{-1}(\R^2)$
of Proposition~\ref{globexist}. A considerable part of the analysis will be
devoted to the continuous dependence with respect to time and the initial data.
For $s\in\R$, we denote by $H^{s,0}(\R^2)$ the anisotropic Sobolev spaces
equipped with the norm
$$
\|u\|_{H^{s,0}(\R^2)}=\|J^s_{x}u\|_{L^2(\R^2)}
$$
For an integer $s\geq 0$ an equivalent norm in $H^{s,0}(\R^2)$ is given by 
$$
\|u\|_{L^2(\R^2)}+\|\partial_x^{s}u\|_{L^2(\R^2)}\,.
$$
We have the following well-posedness result.
\begin{prop}\label{11}
Let $\phi\in H^{2,0}(\R^2)$. Then there exists a unique
  positive $T$ depending only on $\| \phi\|_{H^{2,0}}$ and a unique solution
$u$ of  (\ref{LineKP}) with initial data $\phi$ on
the time interval $[0,T]$ satisfying
\begin{equation}\label{class}
u\in C([0,T];H^{2,0}(\R^2)),\quad  u_x \in L^1_T L^\infty_{xy}\, . 
\end{equation}
Furthermore, the map $ \phi \mapsto u_\phi $ is continuous from $H^{2,0}(\R^2)$ to $C([0,T];H^{2,0}(\R^2))$.
\end{prop}
\subsection{Uniqueness} 
The uniqueness follows from the next lemma.
\begin{lemma}\label{unique}
Let $u, v$ be two solutions of \re{LineKP} in the class defined by (\ref{class}).
Then
$$
\|u-v\|_{L^\infty_T L^2_{xy}}  \leq 
\exp\Bigl( C\|u_x\|_{L^1_TL^\infty_{xy}} + 
C\|v_x\|_{L^1_T L^\infty_{xy}}+CT \|\psi_x\|_{L^\infty} \Bigr)
 \|u(0)-v(0)\|_{L^2_{xy}}\,\, .
$$
\end{lemma}
\begin{proof}
We set 
$$ 
w:=u-v\, . 
$$
Further we define
$$
u^\ep:=\vp\ast u,\quad  v^\ep:=\vp\ast v,\quad w^\ep :=\vp\ast w\, ,
$$
where $ \vp $ is defined by (\ref{b2}). Then $ w^\ep $ solves the equation
\begin{multline*} 
w_t^{\ep}+w_{xxx}^{\ep}-\partial_x^{-1}w_{yy}^{\ep}+\frac{1}{2}
\partial_x\big(w^{\ep}(u+v)\big)+\partial_x(\psi w^{\ep})-
\\
-\partial_{x}\big((\frac{1}{2}(u+v)+\psi)w^{\ep}\big)+\varphi_x^{\ep}\ast
\big((\frac{1}{2}(u+v)+\psi)w\big)=0\,\, .
\end{multline*}
Taking the $L^2$ scalar product of the last equation with $ w^\ep $ gives
\begin{eqnarray*}
 \frac{1}{2}\frac{d}{dt} \|w^\ep\|_{L^2(\R^2)}^2 & = & -\intR
(u_x+v_x) (w^\ep)^2 -\frac{1}{2}\intR \psi_{x} (w^\ep)^2
\\
 & & +\frac{1}{2}
\intR \Bigl(\partial_x [(u+v +2\psi) w^\ep]
  -\vp_x \ast [(u+v+2\psi) w]\Bigr) w^\ep\, .
\end{eqnarray*}
Thus, by Gronwall Lemma,
\begin{eqnarray*}
 \|w^\ep\|_{L^\infty_T L^2_{xy}}^2 & \lesssim & \exp\Bigl( C\|u_x\|_{L^1_T L^\infty_{xy}}
 + C\|v_x\|_{L^1_T L^\infty_{xy}}+CT \|\psi_x\|_{L^\infty} \Bigr)
 \Bigl[\|w^\ep(0)\|^2_{L^2}
\\
 &  & 
+\Bigl\|\partial_x [(u+v +2\psi) w^\ep]
  -\vp_x \ast [(u+v+2\psi) w]\Bigr\|_{L^1_T
 L^2_{xy}} \|w^\ep\|_{L^\infty_T L^2_{xy}}\Bigr]
\end{eqnarray*}
Passing to the limit in $ \ep\rightarrow 0 $, it follows from  Lebesgue theorem and the
assumptions on $u$ and $v$ that the commutator
term in the right-hand side of the last inequality tends to $0$. Thus
\begin{equation}
\|w\|_{L^\infty_T L^2_{xy}}^2  \lesssim \exp\Bigl( C\|u_x\|_{L^1_T
L^\infty_{xy}} + C\|v_x\|_{L^1_T L^\infty_{xy}}+CT \|\psi_x\|_{L^\infty} \Bigr)
 \|w(0)\|^2_{L^2_{xy}} 
 \end{equation}
which completes the proof of Lemma~\ref{unique}.
\end{proof}
\subsection{Existence}
Let $ \phi\in H^{2,0}(\R^2) $. We set
$$
\phi_\ep:=\vp\ast \phi \, ,  
$$ 
where $ \vp $ are defined by (\ref{b2}).
It is clear that $ \phi_\ep\in H^\infty_{-1}\cap Z $ and that
$\phi_\ep\to \phi $ in $ H^{2,0}(\R^2) $.
By Proposition~\ref{globexist}, the emanating solution $ u_\ep $
is global in time and belongs to $ C(\R_+;H^\infty_{-1})$. We will show
that there exists $T=T(\|\phi \|_{H^{2,0}})>0$ such that the sequence 
$\{\partial_x u_{\ep} \} $ is bounded on $ L^1_T L^\infty_{xy} $.
On the level of $H^{2,0}$ regularity we do not control the $L^{\infty}_{xy}$
norm.
Using a Littlewood-Paley decomposition in the $x$ variable and Lemma~\ref{luca}, we can however write
\begin{eqnarray*}
\|u_\ep\|_{L^{1}_T L^\infty_{xy}} & \lesssim & 
\|{\Delta}_1 u_\ep\|_{L^{1}_T L^\infty_{xy}} +\sum_{\lambda\geq 2} \|u_{_\ep} \|_{L^{1}_T L^\infty_{xy}}
\\
& \lesssim &
\|{\Delta}_1 u_\ep\|_{L^{1}_T L^\infty_{xy}} +\sum_{\lambda\geq 2}
\lambda^{-1}\|\partial_x u_{_\ep} \|_{L^{1}_T L^\infty_{xy}}
\\
& \lesssim &
\|{\Delta}_1 u_\ep\|_{L^{1}_T L^\infty_{xy}} + \|\partial_x u_{_\ep}
\|_{L^{1}_T L^\infty_{xy}}\,\, .
\end{eqnarray*}
Next, as in \re{bassefreq}, Duhamel formula and Strichartz estimates of
 Proposition~\ref{Strichartz} yield for $ 0 \le T \le 1 $,
\begin{equation*}
\|{\Delta}_1 u_\ep\|_{L^{1}_T L^\infty_{xy}} \lesssim
\|\phi_\ep\|_{L^2_{xy}} +\|\partial_x u_{\ep}\|_{L^{1}_T
L^\infty_{xy}}\|u_\ep\|_{L^\infty_T
L^2_{xy}}+\|\psi\|_{L^{\infty}}\|u_{\ep}\|_{L^{\infty}_{T}L^{2}_{xy}}\, .
\end{equation*}
Therefore, thanks to the $L^2$ control (\ref{estL2}) of Proposition~\ref{nico1} we obtained
that there exists a constant $C$ depending on bounds on $\psi$ and its
derivatives in $L^{\infty}$ but independent of $\phi$ such that for $0\leq T\leq 1$,
\begin{equation*}
\|u_\ep\|_{L^{1}_T L^\infty_{xy}}\lesssim \|\phi_\ep\|_{L^2_{xy}}
+ 
(1+\|\phi_\ep\|_{L^2_{xy}})\|\partial_x u_{\ep}\|_{L^{1}_T L^\infty_{xy}}\, .
\end{equation*}
Using (\ref{kando1}) and (\ref{kando2}) with $s=2$ we obtain that there exists a
constant $C$ depending on bounds on $\psi$ and its derivatives in $L^{\infty}$ but
independent of $\phi$ such that
\begin{equation}\label{103-bis}
\frac{d}{dt}\|\partial_x^2 u_{\ep}(t,\cdot)\|_{L^2}^{2}\leq 
C\Big(\|\partial_x u_{\ep}(t,\cdot)\|_{L^{\infty}}+C\Big)
\|J^2_x u_{\ep}(t,\cdot)\|_{L^2(\R^2)}^{2}
\end{equation}
and therefore thanks to the $L^2$ bound on $u_{\ep}$ provided by Proposition~\ref{nico1} and the Gronwall lemma, we obtain that
\begin{equation*}
\|J_x^2 u_{\ep}\|_{L^{\infty}_{T}L^2_{xy}}
\lesssim \|J^2_x\phi_{\ep}\|_{L^{2}_{xy}}
\exp(C\|\partial_x u_{\ep}\|_{L^1_{T}L^{\infty}_{xy}}+C\, T )
\end{equation*}
Let us set
$$
f(T):=\|u_\ep\|_{L^{1}_T L^\infty_{xy}}+
\|\partial_x u_{\ep}\|_{L^{1}_T L^\infty_{xy}}\,.
$$
We deduce from the last estimates and (\ref{estim2}) that for $ 0\le T \le 1 $,
\begin{multline}
f(T)\lesssim \, 
(1+\|\phi_\ep\|_{L^2_{xy}})\|\partial_x u_{\ep}\|_{L^{1}_T L^\infty_{xy}}
+
\|\phi_\ep\|_{L^2_{xy}}
\\
\hspace*{-8mm}\lesssim
\|J^{2}_x \phi_\ep \|_{L^2_{xy}}(1+\|\phi_\ep\|_{L^2_{xy}})\exp(Cf(T))\; .
\label{eq10}
\end{multline}
We notice that for $ \alpha>0 $ small enough the continuous  map
$$
g \, :\,  y \longmapsto y-\alpha \exp(Cy) 
$$ 
satisfies $ g(0)<0 $ and that
$ g(y_0)=0 $ for some $ 0<y_0<1 $. Since 
$$ 
T\longmapsto f(T) 
$$ 
is continuous and satisfies $f(0)=0$ and \re{eq10}, we deduce that for
$$ 
\|J^{2}_x \phi_\ep \|_{L^2_{xy}}(1+\|\phi_\ep\|_{L^2_{xy}})
$$ 
small enough, $0\le  f(T)\le 1 $ for all $ T\in[0,1] $.
Therefore, if $\|J^{2}_{x}\phi\|_{L^2}\ll 1$ then for $0\leq T\leq 1$,
\begin{equation}\label{vajno}
\| u_{\ep}\|_{L^1_{T}L^{\infty}_{xy}}+\|\partial_{x} u_{\ep}\|_{L^1_{T}L^{\infty}_{xy}}\leq 1\,.
\end{equation}
We next use a scaling argument to show that a bound of type (\ref{vajno})
holds
for $\|J^{2}_{x}\phi\|_{L^2}$ of an arbitrary size.
Notice that $u(t,x,y)$ is a solution of
$$
(u_t+u_{xxx}+uu_{x}+\partial_{x}(\psi u))_{x}-u_{yy}=0
$$
with $ \phi(x,y)$ as initial data if and only if for every $\beta\in\R$, 
$$ 
u_\beta(t,x,y): =\beta^2 u(\beta^3 t ,\beta x ,\beta^2 y) 
$$
is a solution to  
$$
(\partial_t u_{\beta}+\partial_x^3 u_{\beta}+u_{\beta}\, \partial_x
u_{\beta}+\partial_{x}(\psi_{\beta} u_{\beta}))_{x}-\partial_{y}^{2}u_{\beta}=0,
$$
with initial data
$$
\phi_{\beta}(x,y)=\beta^2\phi(\beta x , \beta^2 y) 
$$
and 
$$
\psi_\beta (t,x, y)=\beta^2 \psi(\beta x - \beta^3 t, \beta^2 y) 
$$ 
instead
of $ \psi (x-t,y) $. One can check that for $0<\beta\leq 1$,
$$
\|J^2_x \phi_{\beta}\|_{L^2_{xy}}\lesssim \beta^{1/2} \|J^2_x
\phi\|_{L^2_{xy}} 
$$
and for $s\in \N$,
$$
\|J^s_x\psi_{\beta}\|_{L^\infty_{xy}}\lesssim \beta^{2} 
\|J^{s}_{x} \psi \|_{L^\infty_{xy}}\,\, .
$$
Let us set
$$
u_{\ep,\beta}(t,x,y):=\beta^{2}u_{\ep}(\beta^3 t,\beta x,\beta^2 y)\, .
$$
In view of the above discussion for 
$$ 
\beta \sim  (1+\|J_x^{2} \phi_\ep\|_{L^2_{xy}})^{-2}
\tendsto{\ep\to 0}  
(1+\|J_x^{2} \phi\|_{L^2_{xy}})^{-2} 
$$ 
one has the bound
$$
\int_0^1
 \|  u_{\ep,\beta} (t) \|_{L^\infty_{xy}}
 +\|\partial_x  u_{\ep,\beta} (t) \|_{L^\infty_{xy}} \, dt
\le 1
$$
which leads to
$$
\int_0^{\beta^3}
\|  u_{\ep}(t) \|_{L^\infty_{xy}} +
 \|\partial_x  u_{\ep}(t) \|_{L^\infty_{xy}}
\, dt \le 1/\beta^{3}
$$
We thus obtain for $ T\sim (1+\|J^{2}_x\phi\|_{L^2_{xy}})^{-6} $ and $ \ep $ small enough,
\begin{equation}
\| u_{\ep}\|_{L^1_T L^\infty_{xy}} +
\|\partial_x u_{\ep}\|_{L^1_T L^\infty_{xy}} \lesssim
(1+\|J^{2}_x\phi\|_{L^2_{xy}})^{6} \label{joji}
\end{equation}
which is the substitute of (\ref{vajno}) in the case of an arbitrary initial data.
Therefore we deduce that $ \{u_\ep \} $ is  bounded in $L^\infty_T H^{2,0}(\R^2)$
and $ \{\partial_x u_{\ep} \} $ is bounded in
$L^{1}_T L^\infty_{xy}$. It then follows from Lemma~\ref{unique} that 
$$ 
\|u_\ep -u_{\ep'} \|_{L^\infty_T L^2_{xy}} \to 0 
$$ 
as $ \ep, \ep' \to 0 $. Hence, there exists
$$
u\in C([0,T];L^2(\R^2))\cap C_w([0,T]; H^{2,0}(\R^2))
$$ 
with $ u_x\in  L^{1}_T L^\infty_{xy}$
such that $u_\ep $ converges to $ u $ in $ C([0,T];L^2(\R^2)) $. It is clear  
that $ u $ satisfies \re{LineKP}  at least in a weak (distributional for
instance) sense. 
\subsection{Continuity with respect to time} \label{63}
Now the continuity of $ u $ with values in $H^{2,0}(\R^2)$ as well as the continuity of the flow-map in $H^{2,0}(\R^2)$ will follow
from the Bona-Smith argument. Note that here we are not in the classical
situation since $ u $ does not belong to $L^\infty(0,T;H^{2+}(\R^2))$. 
Let $u$ be a fixed solution of \re{LineKP} with initial data $\phi\in H^{2,0}(\R^2)$.
Recall that $u_\ep$ is the solution to \re{LineKP} with $\phi_\ep =\vp\ast \phi $ as initial data.
We will show that $ \{u_\ep \} $ is a in fact a Cauchy sequence in $
C([0,T];H^{2,0}(\R^2))$ which will prove that $ u\in
C([0,T];H^{2,0}(\R^2))$ . First by straightforward calculations
in Fourier space, one can show that for 
$ \phi \in H^{2,0}(\R^2) $, $ 0<\ep<1 $ and $ r\geq 0 $,
\begin{equation}\label{approx1}
\|\vp\ast \phi  \|_{H^{2+r,0}(\R^2)} \lesssim \ep^{-r} \|\phi\|_{H^{2,0}(\R^2)}
\end{equation}
and
\begin{equation}\label{approx2}
\|\vp\ast \phi -\phi  \|_{H^{2-r,0}(\R^2)} \lesssim o\big(\ep^{r}\big) \|\phi\|_{H^{2,0}(\R^2)} 
\end{equation}
as $\ep\rightarrow 0$.
For $0< \ep_2 <\ep_1 <1 $, we set 
$$ 
w:=u_{\ep_1}-u_{\ep_2} \, .
$$
It follows from the estimates of the previous subsection applied to
$u_{\ep_1}$ and $u_{\ep_2}$ that
$$
\|w\|_{L^{1}_{T}L^{\infty}_{xy}}+\|w_x\|_{L^{1}_{T}L^{\infty}_{xy}}\leq C
$$
and that for any $s\in\N$, $s\geq 2$,
$$
\|J^{s}_{x}u_{\ep}\|_{L^{\infty}_{T}L^{2}_{xy}}
\lesssim
\|J^{s}_{x}u_{\ep}(0)\|_{L^{2}_{xy}}
\lesssim \ep^{2-s}
$$
The issue is to show that $\|J^{2}_{x}w\|_{L^{\infty}_{T}L^{2}_{xy}}$ is not
only bounded but that it tends to zero as $\ep_1$ tends to zero.
Using Lemma~\ref{unique}, we obtain that for $ 0\leq T\leq 1 $
\begin{equation}\label{w1}
\|w\|_{L^\infty_T L^2_{xy}} \lesssim  \, \|w(0)\|_{L^2_{xy}} 
\lesssim \, o\big(\ep_1^{2}\big)\, .
\end{equation}
Next, we observe that $w$ solves the equation
\begin{equation}\label{W}
w_t+w_{xxx}-\partial_{x}^{-1}w_{yy}+\partial_{x}(\psi
w)-ww_{x}+u_{\ep_1}\partial_{x}w+w \partial_{x}u_{\ep_1}=0\,.
\end{equation}
Note that since $\ep_2 <\ep_1$, we privilege $u_{\ep_1}$ to $u_{\ep_2}$ when
writing the equation (\ref{W}). 
We next apply $\partial_x^2$ to (\ref{W}) and multiply it with $\partial_x^2w$.
Using (\ref{kando1}) and (\ref{kando2}) and  Gronwall lemma, we obtain that
\begin{multline}
\label{w2}
\|J^2_{x} w\|_{L^\infty_T L^2_{xy}}^2  \lesssim  
\exp \Bigl(
C\|\partial_{x}u_{\ep_1}\|_{L^1_T L^\infty_{xy}}+
C\|\partial_{x}w\|_{L^1_T L^\infty_{xy}}+
CT\|\psi\|_{W^{3,\infty}}\Bigr)
\\
\Bigl[ \|J^2_{x} w(0) \|_{L^2}^2  +
\|w_x \|_{L^1_T L^\infty_{xy}}^2 \|J^{2}_x u_{\ep_1}\|_{L^\infty_T L^2_{xy}}^2+\|w \|_{L^1_T L^\infty_{xy}}^2
 \|J^{3}_x u_{\ep_1}\|_{L^\infty_T L^2_{xy}}^2\Bigr]\, .
\end{multline}
On the other hand, according to (\ref{W}), (\ref{estim1}), (\ref{w1}), Lemma 7 and  (\ref{joji}) applied to $ u_{\ep_1} $ and $ u_{\ep_2} $ , one infers that for $ 0\le T \le 1 $,
\begin{eqnarray}
\|w\|_{L^{1}_T L^\infty_{xy}} & \lesssim &
\Bigl(1+\|u_{\ep_1}+u_{\ep_2}\|_{L^1_T L^\infty_{xy}}
+T\|\psi\|_{W^{1,\infty}}\Bigr) \|J_x^{1/2+}
 w  \|_{L^\infty_T
L^2_{xy}} \nonumber \\
&  & +\Bigl( 1+\|D_x^{\frac{1}{2}+}(u_{\ep_1}+u_{\ep_2})\|_{L^1_T
L^\infty_{xy}} \Bigr)
 \| w \|_{L^\infty_T L^2_{xy}} \nonumber \\
  & \lesssim & \ep_1^{\frac{3}{2}-}\label{w4} \quad
\end{eqnarray}
where in the last step we used that, by \re{w1} and interpolation argument,
for
$$
0<\alpha<2
$$
one has the bound
\begin{equation}
\|J^{\alpha}_x w\|_{L^\infty_T L^2_{xy}} \lesssim \|J^{2}_x
w\|_{L^\infty_T L^2_{xy}}^{\alpha/2}\|w\|_{L^\infty_T
L^2_{xy}}^{1-\alpha/2}\lesssim \ep_1^{2-\alpha}\quad .
\label{dec}
\end{equation}
Observing that $ w_x $ solves the equation
$$
(\partial_t+\partial_{x}^{3}-\partial_{x}^{-1}\partial_{y}^{2})w_{x}+\partial_{x}^{2}(\psi
w)+\frac{1}{2}\partial_{x}^{2}(w(u_{\ep_1}+u_{\ep_2}))=0\, ,
$$
another use of  \re{estim1} with \re{dec} and \re{w4} in hands yields
\begin{eqnarray}
\|w_x\|_{L^{1}_T L^\infty_{xy}} & \lesssim &
\Bigl(1+\|u_{\ep_1}+u_{\ep_2}\|_{L^1_T L^\infty_{xy}}
+T\|\psi\|_{W^{2,\infty}}\Bigr) \|J_x^{3/2+}
 w  \|_{L^\infty_T
L^2_{xy}} \nonumber \\
&  & +\|J_x^{\frac{3}{2}+} (u_{\ep_1}+u_{\ep_2})\|_{L^{\infty}_T
L^2_{xy}}
 \| w \|_{L^1_T L^{\infty}_{xy}} \nonumber \\
 & \lesssim & \ep_1^{\frac{1}{2}-}
  +\ep_1^{\frac{3}{2}-}
\nonumber \\
   &  \lesssim & \ep_1^{\frac{1}{2}-} \quad .
  \label{w3}
\end{eqnarray}
Since, according to (\ref{approx1}), $ \|J^{3}_x u_{\ep_1}\|_{L^\infty_T L^2_{xy}}
\lesssim \ep_1^{-1} $, we deduce from \re{w2} and \re{w3} that
for $ 0<T\le 1 $,
\begin{equation*}
\|J^2_{x} w\|_{L^\infty_T L^2_{xy}}  \lesssim \Bigl[
 \|J^2_{x}(\phi_{\ep_1}-\phi_{\ep_2}) \|_{L^2_{xy}}+\ep_1^{0+}
 \Bigr]
\end{equation*}
which proves that $ \{u_{\ep} \} $ is a Cauchy sequence in $
C([0,T];H^{2,0}(\R^2))$. 
\subsection{Continuity of the flow map} \label{64}
Let $ \{\phi_n\}\subset
H^{2,0}(\R^2) $ be such that 
$$ 
\phi_n\longrightarrow \phi 
$$
in $ H^{2,0}(\R^2) $.
We want to prove that the emanating solution $ u^n $ tends to $ u
$ in $ C([0,T];H^{2,0}) $. By the triangle inequality,
$$
\|u-u^n\|_{L^\infty_T H^{2,0}} \le
\|u-u_\ep\|_{L^\infty_T H^{2,0}}+ 
\|u_\ep-u^n_\ep\|_{L^\infty_T H^{2,0}}+
\|u^n_\ep-u^n\|_{L^\infty_T H^{2,0}} \quad  .
$$
Using the argument of the previous subsection, we can obtain the bound
\begin{equation}\label{kak1}
\|u-u_\ep\|_{L^\infty_T H^{2,0}}\lesssim \|u(0)-u_{\ep}(0)\|_{H^{2,0}}+ o(1)=o(1),
\end{equation}
where $o(1)\rightarrow 0$ as $\ep\rightarrow 0$. In the bound for $u-u_{\ep}$
we privilege $u_{\ep}$ to $u$ exactly as we privileged $u_{\ep_1}$ to
$u_{\ep_2}$ in the previous subsection.
Next,  privileging $u^{n}_{\ep}$ to $u_n$, and using that $u^n_{\ep}(0)$
converges to $u^n (0)$ as $\ep\rightarrow 0$, uniformly in $n$, we infer the bound
\begin{equation}\label{kak2}
\|u^n_\ep-u^n\|_{L^\infty_T H^{2,0}} \lesssim \|u^n_{\ep}(0)-u^n (0)\|_{H^{2,0}}+o(1)
= o(1),
\end{equation}
where $o(1)\rightarrow 0$ as $\ep\rightarrow 0$.
Finally, no matter which of the solutions is privileged, we get the bound
\begin{eqnarray}
\nonumber
\|u_\ep-u^n_\ep\|_{L^\infty_T H^{2,0}} & \lesssim &
\|u_\ep(0)-u^n_\ep(0)\|_{H^{2,0}}+o(1)
\\
\label{kak3}
& = & \|\varphi_{\ep}\ast(u(0)-u^{n}(0))\|_{H^{2,0}}+o(1)
\\
\nonumber
& \lesssim & \|\phi-\phi_n\|_{H^{2,0}}+o(1)
\end{eqnarray}
where again $o(1)\rightarrow 0$ as $\ep\rightarrow 0$.
Collecting (\ref{kak1}), (\ref{kak2}) and (\ref{kak3}) ends
the proof of the continuity of the flow map.
Thus the proof of Proposition~\ref{11} is now completed.
\section{End of the proof of Theorem~\ref{thm1} : continuous dependence in $Z$}
In this section, we complete the proof of Theorem~\ref{thm1}.
The existence and uniqueness are obtained as in Proposition~\ref{11}.
Indeed, let us notice that $ Z\hookrightarrow H^{2,0} $ and thus thanks to Proposition~\ref{BorneZ}, the
approximated solutions $u_\ep$ satisfy 
$$ 
\|u_\ep(t)\|_{Z} \le g(t) 
$$ 
 which ensures that  $\|u\|_{H^{2,0}}$ can not go to infinity in finite time.
According to Proposition~\ref{11} , this implies that the solution $u$ is global in time and satisfies
$$ 
u\in L^\infty(0,T;Z) , \quad u_x\in L^1(0,T;L^\infty_{xy}), \quad \forall\, T>0
$$
It remains to prove the continuity of $ t\mapsto u(t) $ with values in
$ Z $ and the continuity of the flow-map. By Bona-Smith type arguments
 as in  the preceding section, one can see that $
t\mapsto u(t) $ as well as the flow-map is continuous with values in $
X $ where
$$ 
X=\{v\in S', \|v\|_{H^{2,0}}+\|\pa u_y\|_{L^2}+\|u_y\|_{L^2} <\infty\, \}
$$
It remains to prove the continuity of $ t\mapsto \partial_x^{-2}
u_{yy} $ with values in $ L^2(\R^2) $. 
This continuity property is more delicate since we can not apply
$\partial_{x}^{-2}$ to the equation satisfied by the approximate solution.
Thus a different argument is needed.
Recall  that 
$$
F^{\psi}(u)=\frac{5}{6} \intR |\partial_{x}^{-2} u_{yy}|^2-\frac{5}{6} \intR
(u^2+2\psi u) \, \partial_{x}^{-2} u_{yy}+r(u)
$$
where $r$ is continuous on $ X $.
\\

We start by proving  that $ t\mapsto F^{\psi}(u(t)) $ is continuous on $
[0,T] $. As in the preceding subsection we denote by $ u_\ep $
the solution emanating from $ \phi_\ep=\vp\ast \phi $. First it is
easy to prove that for any $v\in L^2(\R^2) $ and any $ \alpha> 0$
fixed, the sequence of functions $ t\mapsto (\partial_{x}^{-2}
u_{yy}^\ep(t) , \varphi_\alpha \ast v)_{L^2} $ is equi-continuous
on $ [0,T] $ and thus $
\partial_{x}^{-2} u_{yy}^\ep\to
\partial_{x}^{-2} u_{yy} $ in $ C_w([0,T];L^2(\R^2)) $. Hence, for
any $ t\in [0,T]
 $
$$
\intR (\tu(t))^2 \partial_{x}^{-2} u_{yy}^\ep(t) \to \intR u^2(t)
\partial_{x}^{-2} u_{yy}(t)
$$
$$
\intR |\partial_{x}^{-2} u_{yy}(t)|^2 \le \liminf_{\ep\rightarrow 0} \intR
|\partial_{x}^{-2} u_{yy}^\ep(t) |^2
$$
Thus passing in the limit in $ \ep $ in \re{equF}, we obtain
\begin{eqnarray*}
F^{\psi}(u(t))&  \le & \liminf_{\ep\rightarrow 0} F^{\psi}(u^\ep(t))  
\\
&=& F^{\psi}(\phi)
-\frac{5}{3}\int_0^t \intR \psi_x
(\partial^{-2}_x u_{yy}) \, u_{xx} 
\\
& &
+\frac{5}{6} \int_0^t \intR \psi (\partial_x^{-2} u_{yy})  ((u)^2+2 \psi u)_x
\\
& &
+\frac{5}{3} \int_0^t\intR \psi u (\partial_x^{-2} u_{yy})
+\frac{5}{3} \int_0^t\intR(\psi u)_x u (\partial_{x}^{-2}u_{yy})
\\
& &
+\frac{5}{3}\int_0^t \intR \psi_y(\pa u_y)( \partial_x^{-2} u_{yy})
  +\int_0^t G(u)\,.
\end{eqnarray*}
Taking $ u(t) $ as initial data and reversing time, we get the
reverse inequality and thus
\begin{multline}
\label{equa2F}
F^{\psi}(u(t))  =   F^{\psi}(\phi)-\frac{5}{3}\int_0^t \intR \psi_x(\partial^{-2}_x u_{yy}) \, u_{xx}
\\ 
+\frac{5}{6} \int_0^t \intR \psi (\partial_x^{-2} u_{yy})  ((u)^2+2 \psi u)_x+\frac{5}{3} \int_0^t\intR \psi u (\partial_x^{-2} u_{yy})
\\
+\frac{5}{3} \int_0^t\intR(\psi u)_x u (\partial_{x}^{-2}u_{yy}) +\frac{5}{3}\int_0^t \intR \psi_y(\pa u_y) (\partial_x^{-2} u_{yy}) +\int_0^t G(u) 
\end{multline}
It follows that $ t\mapsto F^{\psi}(u(t)) $ is continuous on  $ [0,T] $.
 Now, since
$\partial_{x}^{-2} u_{yy} \in C_w([0,T];L^2(\R^2)) $,
$$
t\mapsto \intR u^2 \partial_x^{-2} u_{yy} \in C([0,T]) \quad ,
$$
and thus the continuity of $ t\to F^{\psi}(u(t)) $ forces $t\to \intR
|\partial_x^{-2} u_{yy}|^2 $ to belong to $ C([0,T]) $. It
follows that
$$
t\mapsto  \partial_x^{-2} u_{yy} \in C([0,T];L^2(\R^2))) \quad ,
$$
which proves that $ u(t) $ describes a continuous curve in $ Z $.\\
Now, let $ \{\phi_n\} \subset Z $ such that $ \phi_n\to \phi $ in
$ Z $. Since $ F^{\psi} $ is continuous on $ Z $,
  $ F^{\psi}(\phi_n) \to F^{\psi}(\phi) $. Moreover, using that the emanating solutions $ u_n $ converges to $ u $ in $ C([0,T];X)$
  and that $ \partial_x^{-2} u^n_{yy} $ converges to $\partial_x^{-2} u_{yy} $ in $ C_w([0,T];L^2(\R^2)) $,
   we infer from \re{equa2F} that
  $$
  F^{\psi}(u^n(t))-F^{\psi}(\phi_n)\to F^{\psi}(u(t))-F^{\psi}(\phi), \quad \forall t\in [0,T]
$$
This convergence clearly forces
$$
\intR |\partial_x^{-2} u^n_{yy}(t) |^2 \to \intR |\partial_x^{-2}
u_{yy}(t) |^2, \quad \forall t\in [0,T]
$$
which permits to conclude that $ u^n \to u $ in $ C([0,T];Z) $.
\qed
\\

We end this section by an important remark. Notice that our result requires a control
on $\|\nabla_{x,y}\, u \|_{L^1_T L^{\infty}_{xy}}$ in order to have the basic {\bf global}
well-posedness theorem of Section~4. But once this global in time result is
established, further improvements of the local well-posedness theory in the
spaces $H^{s,0}(\R^2)$, or in the spaces considered in \cite{K}, only requires a control on $\|u_x\|_{L^1_T L^{\infty}_{xy}}$
in terms of $\|J^s_x u\|_{L^{\infty}_{T}L^2_{xy}}$.
%
\section{Well-posedness of the KP-I equation in $H^{s,0}(\R^2)$, $s>3/2$}
The aim of this section is to extend Kenig's local well-posedness result by showing that KP-I equation
is locally well-posed for initial data in the space
$$
H^{s,0}(\R^2) :=\{u\in L^2(\R^2)\, :\,  |D_x|^s u\in L^2(\R^2) \}\, \mbox{ with } s>3/2 \, ,
$$
that is no $y$ derivative is needed. Consider thus the KP-I equation
\begin{equation}\label{KP1}
(u_t+u_{xxx}+uu_x)_{x}-u_{yy}=0
\end{equation}
with initial data
\begin{equation}\label{dat}
u(0,x,y)=\phi\in H^{s,0}(\R^2)\, .
\end{equation}
Thanks to the estimates established in Section~4, we have the following modest
extension of Kenig's result \cite{K}.
\begin{Theorem}
The Cauchy problem (\ref{KP1})-(\ref{dat}) is locally well-posed in $ H^{s,0}(\R^2) $ for $s>3/2$.
\end{Theorem}
Let $u$ be a $H^{\infty}_{-1}(\R^2)$ solution to the KP-I equation. Then thanks
to Proposition~\ref{carlos}, we obtain that for $T\leq 1$,
\begin{equation}\label{I}
\|u\|_{L^1_{T}L^{\infty}_{xy}}+\|u_x\|_{L^1_{T}L^{\infty}_{xy}}
\leq
C_{\ep}
\Big(\|J_x^{\frac{3}{2}+\ep}u\|_{L^{\infty}_{T}L^{2}_{xy}}+
\|u\|_{L^1_{T}L^{\infty}_{xy}}\|J_x^{\frac{3}{2}+\ep}u\|_{L^{\infty}_{T}L^{2}_{xy}}\Big)\, .
\end{equation}
Applying $J^s_x$ to (\ref{KP1}) and multiplying it with $J^s_x u$ gives after
applying the Kato-Ponce commutator estimate in $x$,
$$
\frac{d}{dt}\|J^s_x u(t,\cdot)\|_{L^2}^{2}\lesssim
\|u_x(t,\cdot)\|_{L^{\infty}}\|J^s_x u(t,\cdot)\|_{L^2}^{2}
$$
and the Gronwall lemma gives the bound
\begin{equation}\label{II}
\|J^s_x u\|_{L^{\infty}_{T}L^2_{xy}} \leq
\|J^s_x \phi\|_{L^2}\exp(C\|u_x\|_{L^1_{T}L^{\infty}_{xy}})\,.
\end{equation}
Bounds (\ref{I}) and (\ref{II}) enable one to perform a compactness argument
as we did in the proof of Proposition~\ref{compact} which shows that the flow map
of (\ref{KP1})-(\ref{dat}) can be extended to a map on $H^{s,0}(\R^2)$, $s>3/2$ with a
life span depending only on a bound on $\|J^s_x \phi\|_{L^2}$. The continuity of the trajectory in $ H^{s,0}(\R^2) $
 as well as the continuity of the flow-map can be derived as in Sections \ref{63}-\ref{64}


\begin{thebibliography}{99}
\bibitem{APS} J.C. Alexander, R.L. Pego, R.L. Sachs, {\it On the transverse instability of solitary waves 
in the Kadomtsev-Petviashvili equation}, Phys. Lett. A 226 (1997), 187-192.
%
\bibitem{BIN} O. Besov, V. Ilin, S. Nikolski, {\it Integral representation of functions and embedding theorems}, 
J. Wiley, 1978.
%
\bibitem{BS} J.L. Bona, R. Smith, {\it
The initial-value problem for the Korteweg-de Vries equation},
Philos. Trans. Roy. Soc. London Ser. A 278 (1975), 1287, 555--601.
%
\bibitem {Bo} J. Bourgain, {\it
On the Cauchy problem for the Kadomtsev-Petviashvili equation},
GAFA 3 (1993), 315-341.
%
\bibitem{BGT}
N. Burq, P. G\'erard, N. Tzvetkov,
{\it Strichartz inequalities and the nonlinear {S}chr\"odinger equation on compact manifolds},       
Amer. J. Math. 126 (2004) 569-605.
%
\bibitem{Fre} N.C. Freeman,
{\it Soliton interaction in two dimensions},
Advances in Applied Mathematics, 20 (1980), 1-37.
%
\bibitem {FS} A.S. Fokas, L.Y. Sung, {\it
On the solvability of $N$-waves, Davey-Stewartson and
Kadomtsev-Petviashvili equations}, Inverse Problems 8 (1992), 673-708.
%
\bibitem {FP} A.S. Fokas, A.K. Pogorobkov, {\it
Inverse scattering transform for the KP-I equation on the background of a one
line soliton}, Nonlinearity 16 (2003), 771-783.
%
\bibitem{Gallo} C. Gallo,
{\it Korteweg-de Vries and Benjamin-Ono equations on Zhidkov spaces},
Adv. Differential Equations  10  (2005), 277--308.
%
\bibitem{GHS} M.D. Groves, M. Haragus, S.M. Sun, {\it A dimension breaking
phenomenon in the theory of steady gravity-capillary water waves},
 Phil. Trans. Roy. Soc. Lond. A360 (2002), 2337-2358.
%
\bibitem{HP} M. Haragus, R.L. Pego, {\it Travelling waves of the KP equations
with transverse modulations}, C.R. Acad. Sci. Paris 328 (1999), 227-232.
%
\bibitem{IK} A. Ionescu, C. Kenig,
{\it Local and global well-posedness of periodic KP-I equations},
Preprint 2005.
%
\bibitem {IN} R.J. I\'orio Jr., W.V.L.Nunes, {\it
On equations of KP-type}, Proc. Roy. Soc. (1998), 725-743.
%
\bibitem{IM} P. Isaza, J. Mejia, {\it
Local and global Cauchy problems for the Kadomtsev-Petviashvili (KP-II)
equation in Sobolev spaces of negative indices},
Comm. PDE 26 (2001), 1027-1057.
%
\bibitem{JKato} J. Kato
{\it  Existence and uniqueness of the solution to the modified Schr\"odinger map},
Math. Res. Lett.  12  (2005), 171--186.
%
\bibitem {KP} T. Kato, G. Ponce, {\it
Commutator estimates and the Euler and Navier-Stokes equations},
 Comm. Pure Appl. Math. 41 (1988), 891-907.
\bibitem{K} C. Kenig, {\it On the local and global well-posedness for the KP-I
equation}, Annales IHP Analyse Non lin\'eaire 21 (2004) 827-838.
%
\bibitem{KK} C. Kenig, K. Koenig,
{\it On the local well-posedness of the  Benjamin-Ono and modified  Benjamin-Ono equations},
Math. Res. Letters 10 (2003) 879-895.
%
\bibitem{KN} C. Kenig, A. Nahmud,
{\it The Cauchy problem for the Hyperbolic-Elliptic Ishimori system and Schr\"odinger maps},
Preprint 2004.
%
\bibitem {KPV} C.E.~Kenig, G.~Ponce and L.~Vega,
{\it Oscillatory integrals and regularity of dispersive equations}, Indiana Univ. Math. J. 40 (1991), no. 1, p. 33--69.
%
\bibitem{KPV2}{\sc C.E.~Kenig, G.~Ponce and , L.~Vega},
{ \it Well-posedness and scattering results for the generalized
Korteweg-de Vries equation via contraction principle.} Comm. and
Pure Pure and Appl. Math.  46 (1993), 527-620.
%
 \bibitem{KT1}{ H. Koch and N. Tzvetkov}
{\it Local well-posedness of the Benjamin-Ono equation in $ H^s(\R) $},
 I.M.R.N.  26, (2003) 1449-1464.
%
\bibitem{Lions} J.-L. Lions {\it Quelques m\'ethodes de r\'esolution
des \'equations aux d\'eriv\'ees partielles non lin\'eaires},
Dunod, Paris, 1969.
%
\bibitem{mol} L.~Molinet,
{\it On the asymptotic behavior of solutions to the (generalized)
  Kadomtsev-Petvishvili-Burgers equation}, 
J. Diff. Eq. 152 (1999) 30-74.
%
\bibitem{MR}  L.~Molinet and F.~Ribaud,
{\it The global Cauchy probem in Bourgain's type spaces for a dispersive dissipative semilinear equation},
 SIAM J. Math. Anal. 33, (2002), pp. 1269-1296.
%
\bibitem{MST} L. Molinet, J.C. Saut and N. Tzvetkov, {\it
  Global well-posedness for the KP-I equation},
Math. Annalen 324, (2002), pp. 255-275. Correction : Math. Ann. 328 (2004), 707--710.
%
\bibitem{MST-bis} L. Molinet, J.C. Saut, N. Tzvetkov, {\it
Well-posedness and ill-posedness for the Kadomtsev-Petviashvili-I equation},
Duke Math. J. 115(2), (2002), 353-384.
%
\bibitem {Sa} J. C. Saut, {\it Remarks on the generalized
Kadomtsev-Petviashvili equations}, Indiana Univ. Math. J., 42 (1993),
1017-1029.
%
\bibitem{Sch} M. Schwarz, Jr. {\it
Periodic solutions of Kadomtsev-Petviashvili},
Adv. Math. 66 (1987) 217-233.
%
\bibitem{TM} M. Tajeri, Y. Murakami, {\it The periodic soliton 
resonance : solutions of the Kadomtsev-Petviashvili equation with positive
dispersion},
Phys. Lett. A 143 (1990) 217-220.
%
\bibitem {Tak} H. Takaoka, {\it Global well-posedness for the
Kadomtsev-Petviashvili II equation},
Discrete Contin. Dynam. Systems, 6 (2000), 483--499.
%
\bibitem{TT} H. Takaoka, N. Tzvetkov, {\it On the local regularity
of Kadomtsev-Petviashvili-II equation}, IMRN 8 (2001), 77-114.
%
\bibitem{Tom} M. M. Tom, {\it On a generalized
Kadomtsev-Petviashvili equation}, Contemporary Mathematics AMS, 200
(1996), 193-210.
%
\bibitem{Tz} N. Tzvetkov, {\it Global low regularity solutions
for Kadomtsev-Petviashvili equation},
Diff. Int. Eq. 13 (2000), 1289-1320.
%
\bibitem{W} M.V. Wickerhauser, {\it Inverse scattering for the heat equation
and evolutions in $(2+1)$ variables}, Comm. Math. Phys. 108 (1987), 67-89.
%
\bibitem{Z} A.A. Zaitsev, {\it Formation of stationary waves by superposition of solitons},
Sov. Phys. Dokl. 28 (9) (1983), 720-722.
%
\bibitem{ZS} V. Zakharov, E. Schulman, {\it
Degenerative dispersion laws, motion invariants and
kinetic equations}, Physica D 1 (1980), 192-202.
%
%
\bibitem{Zh} X. Zhou, {\it Inverse scattering transform for the time dependent
Schr\"odinger equation with applications to the KP-I equation},
Comm. Math. Phys. 128 (1990) 551-564.



\end{thebibliography}
\end{document}